\documentclass[12pt]{article}

\usepackage[english]{babel}
\usepackage{amsfonts}
\usepackage{epsfig}
\usepackage{cancel}
\usepackage{color}
\usepackage{amsmath}
\usepackage{amsthm}
\usepackage{amssymb}
\usepackage{mathtools}
\usepackage{tikz-cd}
\usepackage{amsmath}
\usepackage{array}
\usepackage{enumerate}
\usepackage{graphicx} 
\usepackage{xcolor}
\usepackage{soul}
\usepackage{mathrsfs}
\usepackage{amssymb}
\usepackage[normalem]{ulem}
\usepackage{mathtools}
\usepackage{accents}
\usepackage{comment}

\newcommand{\om}{{\underline{ {\rm ord}}\, }}

\DeclareMathOperator{\diag}{diag}

\newcommand{\hd}{h_\delta}

\newcommand{\be}{\begin{equation}}
\newcommand{\beqa}{\begin{eqnarray}}
\newcommand{\bl}{\begin{lem}}
\newcommand{\bt}{\begin{teo}}

\newcommand{\C}{\mathbb{C}}
\newcommand{\ch}{\choose}

\newcommand{\dg}{\Delta_{G}}

\newcommand{\ea}{\end{array}}
\newcommand{\ee}{\end{equation}}
\newcommand{\eeqa}{\end{eqnarray}}
\newcommand{\el}{\end{lem}}
\newcommand{\et}{\end{teo}}

\newcommand{\la}{\langle}

\newcommand{\mc}{\mathcal}

\newcommand{\QQ}{\mathbb{Q}}

\newcommand{\nc}{\normalcolor}

\newcommand{\ord}{{\rm ord \,}}

\newcommand{\prm}{{\rm Pr}_{m} }
\newcommand{\prtm}{{\rm Pr}_{2m-1} }

\newcommand{\Z}{\mathbb{Z}}
\newcommand{\ra}{\rangle}

\newcommand{\vnj}{\mathcal{V}^n_j}

\newcommand{\xpq}{x^{P+Q} } 

\newtheorem{teo}{Theorem}
\newtheorem{lem}[teo]{Lemma}
\newtheorem{pro}[teo]{Proposition} 
\newtheorem{cor}[teo]{Corollary} 

\newtheorem{defin}[teo]{Definition}
\newtheorem*{exstar}{Example}             

\theoremstyle{plain}
\newtheorem{remark}[teo]{Remark}

\providecommand{\keywords}[1]
{
  \text{\textit{Keywords---}} #1
}

\title{
{Poincar\'e-Dulac Normal Forms:\\
Formal and Analytic Aspects}}
\author{ Valery G.~Romanovski$^{1,2,3}$ and Sebastian Walcher$^4$ \\
 $^1${\it Center for Applied Mathematics and Theoretical Physics,}\\
{\it University of Maribor,  Mladinska 3, SI-2000 Maribor, Slovenia  } \\  
$^2${\it Faculty of Electrical Engineering and Computer Science,}\\
    {\it University of Maribor, Koro\v ska cesta 160, SI-2000 Maribor, Slovenia}\\
$^3${\it  Faculty of Natural Science and  Mathematics,}\\
{\it   University of Maribor, Koro\v ska cesta 160, SI-2000 Maribor, Slovenia}\\
$^4${\it Fachgruppe Mathematik, RWTH Aachen, D-52056 Aachen, Germany}
}
\date{}

\begin{document}
 \maketitle  

 \begin{abstract} 
{We discuss various aspects concerning transformations of local analytic, or formal, vector fields to Poincar\'e-Dulac normal form, and the convergence of such transformations. After introducing a new transparent approach to transformations, we review Bruno's approach to formal normalization, as well as convergence results in presence of certain versions of Bruno's Condition A. Retracing the proof steps in Bruno's work, we use a different formalism and variants in the line of arguments, which are easily transferrable to more general scenarios. In particular we show how Bruno's approach naturally extends to an elementary proof of Stolovitch's formal and analytic simultaneous normalization theorems for abelian Lie algebras of vector fields. Furthermore we investigate the role of (meromorphic and formally meromorphic) integrability in normal forms and convergence properties. In particular we establish a bridge to Zung's convergence theorems and show that these are can be complemented by Stolovitch's results. 
}
 \end{abstract}


{
\section{Introduction}
\subsection{Background}
{{The present paper is concerned with analytic vector fields in the neighborhood of a stationary point, and their (analytic or formal) normalization.
Normal forms of such vector fields were first introduced and systematically discussed by Poincar\'e \cite{Poincare} and Dulac \cite{Dulac}. The theory was fully formulated and expanded in the seminal work of Bruno \cite{BruB,Bru}}. A good overwiew is given in Anosov and Arnold \cite{ArnAno}. 

Poincar\'e and Dulac proved, for every analytic or formal vector field, the existence of a formal power series transformation to a formal vector field in normal form.  
However, the Poincar\'e-Dulac theorem does not imply that an analytic system admits a convergent normalization, and indeed there exist analytic systems for which no convergent normalization exists. Thus convergence questions turned out to be a central problem in normal form theory.\footnote{For parallel work about smooth normalization of smooth vector fields see e.g.\ Chen \cite{Chen}. } An early convergence proof goes back to Poincar\'e (it is applicable e.g. to  systems for which all eigenvalues of the linearization at the stationary point have negative real parts). Progress was made by Siegel \cite{Siegel} (who introduced diophantine conditions on the eigenvalues), and Pliss (who generalized Siegel's theorem to vector fields in the presence of resonances, but with linear normal form). 
Bruno's \cite{BruB, Bru} comprehensive investigation of normal forms and normalizing transformations identified
two conditions that jointly guarantee convergence for a given vector field were introduced by Bruno in \cite{BruB, Bru}:
\begin{itemize}
    \item A diophantine condition (Bruno's Condition $\omega$) involving the eigenvalues of the linearization at the stationary point, and generalizing Siegel's condition.
    \item An algebraic  condition on the  shape of the normal form; notably Bruno's Condition A (in various forms), which generalizes the condition by Pliss.
\end{itemize}
As was also shown by Bruno, absence of either condition may preclude the existence of a convergent normalization.

In general, for a given vector field there is neither a unique normal form nor a unique normalizing transformation. A particular choice of transformation will be relevant in convergence proofs, for instance. Bruno \cite{BruB} thus introduced and worked with a class of distinguished transformations.
Systematic use of Lie algebra structures and exponentials of vector fields was made by Broer \cite{Broer}, and Walcher \cite{WalcherTraNF}, to preserve special structures. 

 An efficient approach to normal form transformations, together with a useful algebraic formalism, was recently developed by Petek and Romanovski \cite{PR}. It is based on a general discrete Lie-algebraic structure that contains ordinary vector fields as a particular case. Its specialization to the vector fields underlying Bruno's approach to the convergence problem for normal forms leads to what we refer to below as \emph{Bruno vector fields}. 
Bruno vector fields are considerably more than a convenient notation. Their usefulness stems from the fact that they provide a simple and natural Lie-algebraic setting in which formal normalization, resonance relations, integrability, and convergence can be studied within a common framework

Bruno's Condition A is fundamental but it can only give a partial answer to the question what conditions ensure the existence of a convergent normalizing transformation. The importance of symmetries (more precisely, nontrivial commuting vector fields) was noted in Bruno and Walcher \cite{BrW}; for more see the overview in Cicogna and Walcher \cite{CiW}. A definitive characterization is due to Zung \cite{Zung,Zungrat}, who combines symmetries and first integrals (working with radically different methods).

 From a different perspective Stolovitch \cite{St} (see also \cite{St2}) extended Bruno's conditions and theorems to the simultaneous normalization of commuting vector fields. 
 }
 {{
 \subsection{Overview of results}
In the present manuscript we review and discuss various aspects of formal and convergent normalizing transformations,  adding new insights. We will not always consider the most general scenario, or state the most general convergence theorems, but will be content with focussing on some relevant classes. 
The interplay between the algebraic
structure of normal forms, formal integrability, and convergence is 
the main theme of this paper, which contains new, simplified, proofs of known facts as well as new results.

{
In section \ref{sec:trafo} we review and in part re-prove a formal normalization procedure and theorem due to Bruno, using the Bruno vector field formalism which makes arguments more concise and transparent. Moreover we modify some lines of reasoning to lay the groundwork for Stolovitch's theorems. 
In section \ref{sec:converge} we start with recalling Bruno's diophantine and algebraic conditions for convergence of normalizing transformations, adding the observation that the diophantine condition holds whenever the eigenvalues of the linearization are algebraic. We proceed to retrace Bruno's arguments to prove convergence, restricting attention to Bruno's Simplified Condition $A$, and on Bruno's Condition $A_2$ (see Bruno \cite{Bru,BruB}). But here we consider a wider context, including some notions going back to Stolovitch, and we make the line of arguments more transparent due to the formalism introduced in \cite{PR}.
In section \ref{sec:stolo}
we extend Bruno's approach to give an elementary proof of a formal normalization theorem as well as a convergence theorem, originally due to Stolovitch \cite{St}, for commuting vector fields that satisfy an extension of Bruno's Simplified Condition A. While the results here are known, our reasoning requires little technical preparation and is readily accessible to readers who are not familiar with Stolovitch's machinery. Using the Bruno vector field approach here, we obtain a transparent proof of the convergence
theorem. In section \ref{sec:Bruno1} we extend the convergence proofs to a subcase (Condition $A_1$) of Bruno's Case 1, adding some remarks on invariant manifolds.
Finally, in section \ref{sec:integrable} we consider integrability, in the general setting of formally meromorphic first integrals of (formal or analytic) vector fields, and its role for convergence questions. Along the way we re-prove and extend convergence theorems going back to Llibre et al.\cite{LPW} and Zhang \cite{XZ}, and then add a number of new results. For normal forms of vector fields with the maximal formal meromorphic integrability property we give a complete characterization, and we highlight their role in the context of Zung's \cite{Zung,Zungrat} convergence theorems. We close by exhibiting a different version of Zung's convergence theorem, thus building a bridge between Zung's and Stolovitch's work.
}}}


\section{Normal form transformations}\label{sec:trafo}

{
We consider an  autonomous differential system of the form  
		\be
		\label{Asn}
		\dot x =   F(x), \hspace{3mm} x \in  \mathbb{C}^n,
		\ee
		where $F$ is a formal power series with
		$$
		F(x) =Ax +X(x);
		$$
		thus $F(0)=0$, $DF(0)=A$, and $X(x)$ represents the nonlinear terms.

   One of the central topics of this paper is the formal normalization of formal vector fields \eqref{Asn}, and the analytic normalization of analytic vector fields, at the stationary point $0$.   

Throughout we assume that the matrix $A$ is as in a lower-triangular Jordan normal form, but with entries below the diagonal not necessarily equal to 0 or
1. That is,  we can write 
\be\label{eq:mAsig}
A=
\begin{pmatrix}
\lambda_1 & 0 & 0 & \cdots & 0 \\
\sigma_2  & \lambda_2 & 0 & \cdots & 0 \\
0 & \sigma_3 & \lambda_3 & \ddots & \vdots \\
\vdots & \ddots & \ddots & \ddots & 0 \\
0 & \cdots & 0 & \sigma_n & \lambda_n
\end{pmatrix}
\ee 
where 
$
\lambda_i\in \C \ \ (i=1,\dots,n), \ 
\sigma_i\in \C \ \ (i=2,\dots,n).$

Since the matrix \(A\) is in Jordan normal form, a subdiagonal entry
\(\sigma_i\) may be nonzero only if \(\lambda_i=\lambda_{i-1}\).
Conversely, if \(\lambda_i\neq \lambda_{i-1}\), then necessarily
\(\sigma_i=0\).
In the classical Jordan normal form one further normalizes
\(\sigma_i=1\) inside each Jordan block.  Dropping this normalization
does not change the block structure and yields a representation which is sufficient for our purposes.
}}

\subsection{Notation and notions}

For $Q=(q_1, \dots, q_n)\in \Z^n$ we set 
\begin{equation}
    \|Q\|=|q_1+\dots +q_n| \quad \text{  and  }|Q|=|q_1|+\dots +|q_n|. 
\end{equation}
Following the notation in Bruno \cite{BruB} we let  
  \[{\bf N}_i:=\{Q=(q_1, \dots, q_n)\in \Z^n : q_i\ge -1, \,q_j\ge 0 \text {  if \ } j\ne i     \}  \]
  and  
  \[{\bf N}={\bf N}_1\cup \dots \cup {\bf N}_n . \]

  As usual, we abbreviate 
  \[
  x^Q=x_1^{q_1}\cdots x_n^{q_n}.
  \]
  Given $\theta=(\theta_1,\ldots,\theta_n)^\top\in \mathbb C^n$ we consider the term 

 \be \label{Qtm} \Theta=(x\odot\theta)x^Q=\begin{pmatrix}\theta_1x_1\\ \vdots\\ \theta_nx_n\end{pmatrix} \,x^Q,\ee  
where $\odot$ denotes the Hadamard product.
\begin{defin} \label{def:order}
Given term \eqref{Qtm}, 
we  call $Q$ the {\em degree} of this term, and moreover we call $\|Q\|$ the {\em order} of this term and denote $ord\ \Theta= \| Q\| $.
\end{defin}
\begin{remark}
    The notion of degree used in, and fundamental to, this paper is called multi-degree (among other names) by other authors (see e.g. \cite{MS}). The notion of order is compatible with the Lie algebra structure (to be introduced below).
\end{remark}
  
{
In his 1962 work \cite{Bru62} Bruno proposed to write  analytic or formal vector  fields in the form 
\be \label{B_or}
\dot x_i= x_i f_i(x)= x_i \sum_{Q\in {\bf N}_i} f_{iQ} x^Q,\quad i=1,\dots, n.
\ee

With this notation we will (as was first proposed in \cite{PR})  write  \eqref{Asn} as {\em Bruno vector field} in the form
\be \label{sys}
\dot x = F(x)=\sum_{Q\in{\bf N}} (x\odot F_Q) x^Q,
\ee 
using the Hadamard product, where
\[
F_Q=(F_Q^{(1)},\dots, F_Q^{(n)})^\top\in\mathbb C^n, \quad   F_{\bar 0}=\lambda:=(\lambda_1,\ldots,\lambda_n)^\top.
\]
}

From now on we assume that $\bf N$ is ordered using { the degree  lexicographic 
 order with $x_1>\dots >x_n$, that is for $P=(p_1,\dots, p_n), Q=(q_1,\dots, q_n)\in {\bf N}$ one sets  $P\prec Q$ if the first non-zero  number in the sequence
 $$
 \|Q\|-\|P\|, q_1-p_1, \dots, q_{n-1}-p_{n-1}
 $$
 is positive, 
and that system  \eqref{sys} is   written according to  this  order.}


\begin{remark}
Denote by $e_1,\ldots,e_n$ the standard basis of $\mathbb C^n$. Then for a single vector monomial we have
\begin{equation}\label{notation}
   x_1^{m_1}\cdots x_n^{m_n}e_k=(x\odot e_k) x^Q,\quad Q=(m_1,\ldots, m_k-1,\ldots,m_n)\in {\bf N}_k
\end{equation}
\end{remark}

\begin{defin}
\begin{itemize}
\item  { For $Q=(q_1,\dots,q_n)\in {\bf N}$  we set 
$$ \left<Q,\lambda\right>=\sum q_i\lambda_i.
$$
}
\item The term $  (x\odot F_Q) x^Q $ of \eqref{sys} is called resonant if and only if $\la Q,\lambda\ra =0$.
\item 
			We say that the vector field \eqref{sys} is in normal form up to degree $P$ if
for all $Q\le P$, $F_Q=\bar 0$ when $\left<Q,\lambda\right>\not=0$. 		
\item 				We say that the vector field \eqref{sys} is in (Poincar\'e-Dulac) normal form if
for all $Q\in {\bf N}$, $F_Q=\bar 0$ when $\left<Q,\lambda\right>\not=0$. 	
\end{itemize}
\end{defin}

We denote
\begin{equation}\label{efre}
{M}_\lambda:=\left\{\alpha \in {\bf N}| \  \langle\alpha,\lambda\rangle=0
\right\}.
\end{equation}
Recall that 
\be \label{linp}
A x= (x\odot \lambda ) x^{\bar 0}+  (x\odot (0,\dots, 0, \sigma_n))^\top x_{n-1} x_n^{-1} 
+\dots +   (x\odot (0,\sigma_2,0, \dots, 0))^\top  x_1 x_2^{-1},
\ee
where 
$\lambda=(\lambda_1, \dots, \lambda_n)^\top $.
Thus, if $\sigma_k $ in \eqref{linp} is non-zero, then the term $$ (x\odot (0,\dots, 0, \sigma_k,0,  \dots, 0 ))^\top x_{k-1} x_k^{-1} $$ is resonant.
Denoting this term by $J_k(x)$ we can write \eqref{linp} as 
\be \label{linp7}
A x= (x\odot \lambda ) x^{\bar 0}+ \sum_{k=2}^n J_k(x).
\ee

Therefore,  according to the above definition system \eqref{sys} is in normal form if and only if it is written as 
 $$
\dot x =\sum_{Q\in M_\lambda} (x\odot F_Q) x^Q. 
$$

\begin{defin} Let $\mathcal M$ be the set of all vector fields of the form \eqref{sys}{(note that every formal power series vector field can be expressed in this notation). } Moreover:
\begin{itemize}
    \item For $Q\in{\bf N}$ we let $\mathcal{M}_Q\subset\mathcal  M$ be the vector space of homogeneous vector fields of degree $Q$, 
that is,
$$
\mathcal{M}_Q=\{ (x\odot  \theta) x^Q \ : \ \theta\in \C^n \}.
$$
\item For $\Theta\in \mathcal M$ we say that $m_1\leq \text{ord\, } \Theta < m_2$ if all terms of $\Theta$ are of order greater or equal than $m_1$ and of order $<m_2.$ 
\item We denote by $\mathcal{V}^{m}$ the space of vector fields of order $\leq m$ and by  $\mathcal{V}^n_{j}$ the subspace of vector fields of order $j$, for all $j\leq m$.
\end{itemize}

\end{defin}
		With these notions we will also use a variant notation for \eqref{sys}, viz.
\begin{equation}\label{sysord}
\dot x = F(x)=Ax+\sum_{j\geq 1} F_j(x),
\end{equation}
where $F_j$ denotes the homogeneous term of order $j$ in the expansion, so $F_j\in \mathcal{V}^n_j$ whenever $j\leq m$.


\subsection{Lie algebra structure}
As is well known, the formal vector fields form a Lie algebra with respect to the bracket
\[
\left[R,\,S\right](x)=DS(x)\,R(x)-DR(x)\,S(x).
\]
We note some observations: Every linear vector field $B$ induces the adjoint action
\[
{\rm ad}_B: V\mapsto \left[B,\,V\right];\quad \left[B,\,V\right](x)=DV(x)\,B(x)-B\,V(x)
\]
on vector fields, and an action 
\[
L_B: \psi\mapsto L_B(\psi); \quad L_B(\psi(x)=D\psi(x)\,B(x)
\]
as derivation on functions. We recall some facts about the adjoint action. 
\begin{lem}\label{linaclem}
    Let $B=\diag(\mu_1,\ldots,\mu_n)$ be diagonal, thus $Bx=(\mu_1,\ldots,\mu_n)^\top\odot x=\mu\odot x$. Then:
    \begin{itemize}
        \item The eigenvalues of the linear action ${\rm ad}_B$ on $\mathcal V^m$, the space of all vector fields with order $\leq m$, are the $\left<Q,\,\mu\right>$, with $Q\in\bf N$, $\|Q\|\leq m$, and corresponding ``eigenvectors'' $(x\odot F)x^Q$, $F\in \mathbb C^n$. Moreover,
        $\mathcal V^m$ is the direct sum of the eigenspaces.
        \item If $R_1,\,R_2\in \mathcal V^m$ with $[B,\,R_i]=\delta_iR_i$, then $[B,\,[R_1,\,R_2]]=(\delta_1+\delta_2) \,[R_1,\,R_2]$.
\item If $R\in \mathcal V^m$ with $[B,\,R]=\sigma R$, $\sigma\in\mathbb C$, and $\psi$ is a polynomial with $L_B(\psi)=\rho \psi$, $\rho\in\mathbb C$, then 
\[
[B,\,\psi R]=(\sigma+\rho)\,\psi R.
\]
    \end{itemize}
\end{lem}
\begin{proof}[Sketch of proof]
The second statement follows from the Jacobi identity. For the third, evaluate the Lie bracket.
\end{proof}
\begin{defin}\label{def:eigen} Let $B$ and $\mu$ be as above.
For $\delta\in \mathbb C$ we denote by  $\mathcal V^m_\delta$ the $\delta$-eigenspace of ${\rm ad}_B$, or  $\mathcal V^m_\delta=\{0\}$ if $\delta$ is not an eigenvalue of ${\rm ad}_B|_{{\mathcal V}^m}$.\\
For
\[
R=\sum_Q(x\odot R_Q)\,x^Q\in \mathcal V^m
\]
we set
\[
R_\delta=\sum_{Q:\,\left<Q,\mu\right>=\delta}(x\odot R_Q)\,x^Q,
\]
the component of $R$ in $\mathcal V^m_\delta$.
\end{defin}

We now consider the Lie algebra structure from the perspective of vector fields written in the form \eqref{sys}.
\begin{itemize}

\item 
The following lemma (a particular case of Lemma 2 of \cite{PR}) provides an explicit expression for the Lie bracket of vector monomials. 
 
\begin{lem}\label{lem_tpm}
			If ${ \Theta}=    (x \odot \theta_\mu ) x^\mu   $ and 
			$ { \Phi} =  (x \odot  \phi_\nu )  x^\nu   $, 
			where $\mu, \nu \in {\bf N}$, $\theta_\mu = (  \theta_1 , \dots,  \theta_n  )^\top \in \mathbb C^n$ 
			and   $\phi_\nu = (  \phi_1  ,\dots,  \phi_n  )^\top\in\mathbb C^n$,   then 
			\be\label{tpm} 
			[{ \Theta},{ \Phi}]= 
		(x \odot 	 \alpha_{\mu +\nu})x^{\mu+\nu}  \in \mathcal M, 
			\ee
			with 
			\be \label{alphamunu}
			\alpha_{\mu +\nu}=\la \nu ,\theta_\mu \ra { \phi_\nu}    -    \la  \mu , {\phi_\nu}\ra  \theta_\mu.
			\ee
            {
            In particular one has ${\rm ord}\,([\Theta,\Phi])\geq {\rm ord}(\Theta)+{\rm ord}(\Phi)$.
            }
		\end{lem}


Recall that 
a graded Lie algebra is a Lie algebra that has a decomposition into homogeneous components indexed by a monoid  and the Lie bracket respects this grading.
We grade $\mathcal M$ using $\Z^n$ defining the homogeneous components $\mathcal M_Q$ ($Q\in \Z^n)$ by 
$$
{\mathcal M}_Q= \left\{\begin{matrix}
{\bf 0},  \  {\rm if \ } Q\not \in {\bf N}\\
(x \odot F) x^Q , F\in \C^n   {\rm\  if \ } Q \in {\bf N}
\end{matrix} \right.
$$  
Then any formal series \eqref{sys} is written as a sum of homogeneous components 
and it holds that
\be
\label{Mpq}
[{\mathcal M}_P, {\mathcal M}_Q]\subset {\mathcal M}_{P+Q}.
\ee
\item 	
		By Lemma \ref{lem_tpm} if ${\bf \Theta} \in  \mathcal{M}_P $,
		${ \Phi} \in  \mathcal{M}_Q $ 
		and  $ \text{ad}_{\Theta}$  is the   adjoint operator acting on $ \Phi$
		by
		$$
		(\text{ad}_{\Theta}) \, { \Phi} =[{ \Theta},{ \Phi}],
		$$
		then  $ (\text{ad}_{ \Theta})^i { \Phi}$ is an element of 
		$\mathcal{M}_{Q+iP} $, that is,
		\be\label{ad_i}
		(\text{ad}_{\Theta})^i :  \mathcal{M}_{Q} \nc  \to  \mathcal{M}_{Q+iP}. 
		\ee 
\item 
We introduce an operator for later use. First, we note:
\begin{lem}\label{lem:Du}
For $U ,V \in \mathcal M$ such that 
$$
U=\sum_Q (x \odot   U_{Q}) x^{Q},\quad V=\sum_P (x \odot   V_{P}) x^{P},
$$
one has
\be \label{DU1}
(D U) V=  \sum_{Q,P}(x\odot U_Q\odot V_P) x^{P+Q} +\sum_{Q,P} \la Q, V_P\ra (x \odot U_Q) x^{P+Q}.  
\ee
\end{lem}
\begin{proof}
We recall that  for $f: \C^n \to \C^n$ the derivative of the Hadamard product reads
			$$
			D (x \odot f(x))= {\rm diag}(f)  +{\rm diag}(x) Df(x). 
			$$
Then, 
\begin{equation*}
\begin{array}{rcl}
     D(x \odot (U_Q x^Q))  (x \odot (V_P x^P)) &=&[{\rm diag}(U_Q x^Q)+ {\rm diag}(x) D(U_Q x^Q )  ]  (x \odot (V_P x^P))\\
&=&{\rm diag}(U_Q x^Q) (x \odot (V_P x^P)) \\
 & &+ {\rm diag}(x) D(U_Q x^Q )   (x \odot (V_P x^P)) \\ 
&=&(x\odot U_Q\odot V_P) x^{P+Q} +
 \la Q, V_P\ra (x \odot U_Q) x^{P+Q}.
 \end{array}
\end{equation*}
\end{proof}
\begin{defin}\label{Deldef}
    Given $F(x):=x\odot f(x)$, we set
    \begin{equation}
        \Delta_{F}(x):={\rm diag}(x)\cdot Df(x).
    \end{equation}
\end{defin}
\begin{remark}
 By construction, we have 
 \begin{equation}\label{DvsDelt}
     DF(x)={\rm diag}\, f(x)+ \Delta_{ F}(x).
 \end{equation}
 Moreover, with $G(x):=x\odot g(x)$ one finds
 \begin{equation}\label{LieDel}
  \left[F,\,G\right](x)=   DG(x)F(x)-DF(x)G(x)=\Delta_G(x)F(x)-\Delta_F(x)G(x).
 \end{equation}
\end{remark}

{
For  $U$ and $V$ as in Lemma \ref{lem:Du}
by $ \Delta_U $ we denote the operator which acts on $V\in \mathcal{M}$ by 
\be \label{DelU}
\Delta_UV:=\sum_{Q,P} \la Q, V_P\ra (x \odot U_Q) x^{P+Q}.
\ee

The following statements are then obvious.
 
\begin{lem}\label{lem:9} For  $U,V \in \mathcal{M}$, 

(a) $[U,V]=\Delta_VU-\Delta_UV= {\sum_{P,Q} \la P, U_Q\ra (x \odot V_P) x^{P+Q}.  -\sum_{Q,P} \la Q, V_P\ra (x \odot U_Q) x^{P+Q};}$

(b) {
$
(D U) V=  \sum_{Q,P}(x\odot U_Q\odot V_P) x^{P+Q} +\Delta_U V.
$}
\end{lem}
}
\begin{lem}\label{lem:Deltasq}

For $U ,G \in \mathcal M$ such that 
$$
U=\sum_Q (x \odot   U_{Q}) x^{Q},\quad G=\sum_P (x \odot   G_{P}) x^{P},
$$
 we have 

\be \label{DG2} 
\begin{aligned}
\dg^2(U)=\dg(  \sum_{P,Q} \la P, U_Q \ra (x \odot G_P) x^{ P+Q})= 
\\
\sum_{P_2}\sum_{P_1}\sum_{Q} \la P_1, U_Q \ra \la P_2, G_{P_1}\ra (x \odot G_{P_2}) x^{P_2+P_1+Q}. 
\end{aligned}
\ee
\end{lem}

\end{itemize}

\begin{remark} Vector fields in the form \eqref{sys} were introduced in the works of A.D.\  Bruno, so it is reasonable  to call them Bruno vector fields.  However he wrote them not  in the form \eqref{sys} but as  
$$
\dot x_i= x_i f_i(x)= x_i \sum_{Q\in {\bf N}_i} f_{iQ} x^Q,\quad i=1,\dots, n.
$$
The notation using the Hadamard product and the Lie algebra structure were introduced in  
 \cite{PR}.  
\end{remark}


\subsection{Transformations}\label{subsec:trafos}
It is well known  that  any   formal  change of coordinates of the form 
		\be
		\label{yhy}
		x = H(y) =y+ \  h(y)= y +\sum_{j=1}^\infty  h_j(y), 
		\ee
		with  $ h_{j} \in 
		\mathcal{V}^n_{j}$ for all $j\geq 1$, 
		brings system \eqref{Asn} to a system of a similar form, 
		\be
		\label{linearni}
		\dot y = \widetilde F(y)=A y+ Y(y),
		\ee
		where 
		$$ Y(y)= \sum_{j=1}^\infty \ \widetilde F_j(y), \quad \mbox{ and }  \widetilde F_j(y)\in \vnj \quad \mbox{ for all } j \geq 1.
		$$ 
We will take a closer look at such transformations now.\\
Given $H$, the right hand side of the transformed system is determined by 
\begin{equation}\label{sopremap}
DH(y)\,\widetilde F(y)=F(H(y)).
\end{equation}
The following observations provide some detail about low order terms on both sides of \eqref{sopremap}, and show that they are of relatively simple shape.
\begin{defin}
Given a formal vector field $R=\sum_{j\geq -1} R_j$ of the form \eqref{sys}, with $R_j$ homogeneous of order $j$, and a nonnegative integer $k$, we denote by
\begin{equation}
{\rm Pr}_k(R)=\sum_{j=-1}^k R_j
\end{equation}
its projection to the space $\mathcal {V}^k$.
\end{defin}
\begin{lem}\label{lowlem}
Assume that $m\geq 1$ and that $H$ is given by \eqref{yhy}, with $h$ containing only terms of order $\geq m$. Then for the left and right hand side of \eqref{sopremap} one has
\[
{\rm Pr}_{2m-1}\left(DH(y)\widetilde F(y)\right)={\rm Pr}_{2m-1}\left(\widetilde F(y)+Dh(y) \widetilde F(y)\right),
\]
and 
\[
{\rm Pr}_{2m-1}\left(F(H(y))\right)={\rm Pr}_{2m-1}\left(F(y)+DF(y)\,h(y)\right).
\]
\end{lem}
\begin{proof} As a matter of terminology in this and the following proof, for a scalar polynomial we speak of its conventional degree. For a vector field \eqref{Qtm} we call its conventional order the highest conventional degree of its components (which equals $1+\|Q\|$). \\
For the left hand side there is nothing to show. There remains the right hand side.
\begin{enumerate}[(i)]
\item Let $S(y)$ be homogeneous of order $k$. Then for any $\ell\leq k+1$ one may view the $\ell^{\rm th}$ derivative $D^\ell S(y)$ as a symmetric $\ell$-linear map on $\mathbb C^n\times\cdots\times \mathbb C^n$ (a tensor), with coefficients of conventional degree $k-\ell+1$. (For $\ell>k+1$ this is the zero map.) Given that $h$ contains only terms of order $\geq m$, to compute the order of the terms of $ D^\ell S(y)\,\left(h(y),\ldots,h(y)\right)$ 
we first note that the sum of the conventional orders of terms in $\left(h(y),\ldots,h(y)\right)$ is at least  $m\ell +\ell$. Thus,  the conventional order of the terms of    $ D^\ell S(y)\,\left(h(y),\ldots,h(y)\right)$ is $\geq  m\ell +k + 1$. 
Since, by Definition \ref{def:order},
the order is  one less than the conventional order,  it follows    that
$D^\ell S(y)\,\left(h(y),\ldots,h(y)\right)$ contains only terms of order $m\ell +k$ or higher.

\item Now  $F(x)=Ax+\sum_{j\geq 1} F_j(x)$ with homogeneous $F_j$. With Taylor we have
\[
F_j(y+z)=F_j(y)+DF_j(y)\,z+\frac12 D^2F_j(y)(z,z)+\cdots,
\]
hence
\[
F_j(y+h(y))=F_j(y)+DF_j(y)\,h(y)+\cdots,
\]
where all the summands $D^\ell F_j(y)\,\left(h(y),\ldots,h(y)\right)$, $\ell\geq 2$, contain only terms of orders 
\[
m\ell+j\geq 2m.
\]
\end{enumerate}
\end{proof}
With similar arguments we obtain a further auxiliary result.
\begin{lem}\label{comblem}
Let $m\geq 1$ and $h_1,\ldots,h_q$ of order between $m$ and $2m-1$.  Then for
\[
H^{[k]}(x)= x+h_k(x),\quad 1\leq k\leq q,\text{  and  } H:=H^{[q]}\circ\cdots\circ H^{[1]}
\]
one has
\[
{\rm Pr}_{2m-1}\left(H(x)\right)=x+h_1(x)+\cdots+h_q(x)
\]
\end{lem}
\begin{proof}
        For $j\geq 1$ and any $y$ Taylor expansion yields
        \[
        H^{[j]}(x+z)=H^{[j]}(x) + DH^{[j]}(x)z+\frac12D^2H^{[j]}(x)(z,z)+\ldots,
        \]
        where
        \[
        DH^{[j]}(x)={\rm Id}+Dh_{j}(x),\quad D^\ell H^{[j]}(x)=D^\ell h_{j}(x)\text{  for }\ell\geq 2.
        \]
        Thus all coefficients of $D^\ell H^{[j]}(x)$ have conventional degree $\geq 1+m-\ell$. Hence for any $k\geq m$ the terms in $D^\ell H^{[j]}(x)(h_k(x),\ldots,h_k(x))$ are of conventional order greater than or equal to
        \[
       1+m+(k+1-1)\ell\geq 1+m+m\ell>2m,
        \]
        hence of order $>2m-1$.
        Moreover the only summand in 
        \[
        DH^{[j]}(x)h_k(x)=h_k(x)+\cdots
        \]
        which may contain a term of order $\leq 2m-1$ is equal to $h_k(x)$.
        Therefore by induction we get
        \[
        H^{[j]}\circ\cdots\circ H^{[1]}(x)=x+h_1(x)+\cdots+h_{j}(x)+\cdots
        \]
        for all $j\leq q$, with the dots indicating terms of order $\geq 2m$.
\end{proof}
        

\subsection{Normal form computations}
The following result is  {essentially Theorem 1 of Bruno \cite[p.153]{BruB}. The proof  we give is  based on a different formalism.
 Recall that by our assumption terms of \eqref{sys}  are ordered by degree lexicographic order.}
{   	\begin{teo}\label{th_nf2}
Assume  that the vector field \eqref{sys} is in normal form up to degree  $Q$.
Let $S\in {\bf N}$, $Q< S$,  $F_S\ne \bar 0$ and assume that there is no $T\in {\bf N}$
such that $Q< T < S$ and $F_T\ne \bar 0$; thus the  smallest  term of degree $>Q$  has degree $S$.
	
1) If $\la S,\lambda\ra =0$, then \eqref{sys} is in normal form up to degree
$S$. 

2) If $\la S,\lambda\ra \ne 0$, let
\be \label{hsd2}
h(y)=  \frac 1{  \la S, \lambda \ra  }  (y \odot  F_S )   y^S 
\ee
and
\be \label{Hs2}
			x=H(y):=y+h(y).
			\ee  
Then performing    substitution \eqref{Hs2} we obtain from \eqref{sys} a vector field which is in			normal form up to degree $S$.
		\end{teo}
		\begin{proof}
Let
    $h(y)=   (y \odot  H_S )   y^S. 
$
Then after the transformation $x=y+h(y)$ we obtain from \eqref{sys} the equation 
\be \label{eq:ydbt1}
\begin{aligned}
\dot y= &({\rm Id}+D h(y ))^{-1} F(y+h(y))=
({\rm Id}-D h(y) +\dots  ) (A (y+h(y))+\widehat F(y+h(y))+\dots  )\\
= &  A y +( ( y \odot F_S) y^S -[Ay, h(y)])+\dots ,
\end{aligned}
\ee
where the dots stand for terms of degree greater than  $S$ and $\widehat F$ denotes the nonlinear part of $F$.
We briefly write $[\lambda,h ]$ instead of $[ \lambda \odot y, h(y)]$.
From \eqref{linp7} we have 
$$
[Ay, h(y)]=[\lambda,h]+\sum_{k=2}^n  [J_k(y), h(y)]. 
$$
Taking into account \eqref{tpm} we see  that the terms under the summation sign in the above formula have degree grater than $S$.
Thus the vector field \eqref{eq:ydbt1} does not have a term of degree $S$ if 
\be \label{eq:hs2}
(y \odot F_S) y^S-[\lambda  , h]=0
\ee
Therefore choosing  $h(y) $ as in \eqref{hsd2} we obtain vector field \eqref{eq:ydbt1}  in normal form up to degree $S$.
\end{proof}
}

By induction one obtains the Poincar\'e-Dulac theorem (which we include for the sake of completeness):
\begin{cor}	
For any vector field \eqref{sys} there exists an invertible formal transformation to normal form.
\end{cor}


In general such a transformation is not unique, but there is a unique transformation whose expansion contains no resonant terms.  Bruno \cite{Bru} (Ch.\,III, section 1.5) calls this the distinguished transformation.

We now look more closely at normalizing transformations and the way they are constructed.
The following statement, which extends a result of Pliss \cite{P} to the case of nonlinear normal forms, is due to Bruno \cite{BruB}; \S4, Lemma 3,
but our proof is simpler.

\begin{pro}\label{pro:1n}
Assume that system \eqref{sys} coincides with the normal form up to terms of order $m-1$, 
that is,{
\be \label{nfm}
\dot x=F(x)=\lambda \odot  x +G_*(x)+F_*(x),
\ee
}
where 
$$
{0} \le \text{ord\, } G_*(x)\le m-1, \quad m\le \text{ord\, } F_*(x) 
$$
and $G_*(x) $ contains only resonant terms {of degree greater than $\bar 0$}.
There exists $h\in \mathcal M $ with
\be\label{cond:hm}
m\le \ord h \le 2m -1 
\ee
such that the substitution 
\be\label{xyh}
x=y+h(y)
\ee
 transforms  \eqref{nfm} into a system 
\be \label{nfmt}
\dot y= \widetilde F(y)=G(y)+ \widetilde F_*(y) 
\ee
where
\be \label{GsGt}
G(y)= \lambda \odot   y+  G_*(y)+  \widetilde G(y)
\ee
and
$$
m\le \ord \widetilde G(y)\le 2m-1, \quad 2 m\le \ord \widetilde F_*(y), $$
and  moreover $\widetilde G$ is the sum of the resonant terms in ${\rm Pr}_{2m-1}(F_*)$.
\end{pro}

\begin{proof} 
We order  all $Q$,  $ m\le \|Q\| \le 2m -1$ using the degree lexicographic ordering. We write 
\[
{\rm Pr}_{2m-1}(F_*)=\widetilde G+\widehat F_*.
\]
Let $Q_1$ be the smallest degree that appears in $\widehat F_* $, and 
write 
\be \label{hprmn3}
\widehat F_*(x)=\sum_{i=1}^u (x\odot \widehat F_{*,Q_i}) x^{Q_i},
\ee 
where the sum extends over  all $n$-tuples of $\bf N$ such that $ m\le \|Q_i\| \le 2m -1$ (so some $\widehat F_{*,Q_i}$ in \eqref{hprmn3}  may be equal to zero). By Theorem \ref{th_nf2}, the substitution
$x=H^{[1]}(y):=y +h_1(y),$
where 
$$
h_1(y)=\frac{1}{\la Q_1,\lambda\ra}(y\odot F_{*,Q_1})\,y^{Q_1},
$$
yields a transformed system
 \[
\dot y=\check F(y)=\lambda \odot  y+\check G_*(y)+\check F_*(y),
\]
with $1\leq\ord\check G_*\leq m-1$, $m\leq \ord \check F_*$.
We use Lemma \ref{lowlem} to compare terms of order $<2m$. One sees
\[
{\rm Pr}_{2m-1}\left(DH^{[1]}(y)\check F(y)\right)={\rm Pr}_{2m-1}\left(\lambda \odot  y+\check G(y)+\check F_*(y)+ Dh_1(y)\, (\lambda \odot y)+D h_1(y)\check G(y)\right),
\]
because $Dh_1(y)\,\check F_*(y)$ contains only terms of order $\geq 2m$.
Moreover
\[
{\rm Pr}_{2m-1}\left(F(H^{[1]}(y))\right)={\rm Pr}_{2m-1}\left(\lambda \odot y+G_*(y)+F_*(y)+ \lambda \odot  h_1(y)+ DG_*(y)\,h_1(y)\right),
\]
because $DF_*(y)h_1(y)$ contains only terms of order $\geq 2m$. Comparing terms of like order we find
\[
\check G_*=G_*
\]
and furthermore {
\[
\check F_*=F_*-[h_1,\,\lambda] -\sum_{i=2}^n [h_1, J_i]   -[h_1,\,G_*].
\]
}
By construction of $h_1$, the degree $Q_1$ term of $\check F_*$ is { non-resonant}, hence $\check F$ is in normal form up to degree $Q_1$. Moreover, each of its terms of degree $Q_i>Q_1$ has the same resonant term as the corresponding term in $F_*$, since $[h_1,\,G_*]$  {and  $[h_1, J_i]$ ($i=2,\dots, n$) contain } no resonant terms; see equation \eqref{Mpq}.\\

Iterating the procedure for $Q_2,\ldots,Q_u$ we obtain transformations
$$
H^{[i]}(y)= y+h_i(y)
$$
such that $H=H^{[u]}\circ\cdots\circ H^{[1]}$ transforms 
 \eqref{nfm} into the system \eqref{nfmt}, and $\widetilde G$ is the sum of the resonant terms in ${\rm Pr}_{2m-1}(F_*)$.
Finally, Lemma \ref{comblem} shows that the Taylor expansion of $H$ is of the form
\[
H(y)=y+h_1(y)+\cdots+h_u(y)+\cdots,
\]
where the dots stand for terms of order $\geq 2m$. With $h=h_1+\cdots +h_u$ and again with Lemma \ref{lowlem}, the proof is finished.
\end{proof}

\begin{remark}  The proof of the Proposition is based on the same idea as the proof of Lemma 3  in Bruno \cite[\S 4]{BruB},
but it uses a different formalism and proceeds step-by-step in degrees (rather than Bruno's approach in one leap). This approach will carry over easily to our proof of a convergent normalization result due to Stolovitch (see Section \ref{sec:stolo} below).
\end{remark}

\begin{cor}\label{cor47}
The function  $h$  in Proposition \ref{pro:1n} can be written as 
$$
h=\sum_{\delta\ne 0} h_\delta
$$
where for each eigenvalue $\delta$ of ${\rm ad}_A$, 
$$h_\delta=\sum_{Q: \la Q, \lambda\ra=\delta }  (x\odot H_Q) x^Q
$$  is the unique solution of
\begin{equation}\label{ehds}
\delta\,h_\delta+{\rm Pr_{2m-1}}\left(\left[ G_*, \,h_\delta\right]\right)={\rm Pr}_{2m-1}(F_{*,\delta}) \text{  in  } \mathcal V^{2m-1}.
\end{equation}
\end{cor}
\begin{proof} 
From the Proposition, there is a substitution $x=H(y)=y+h(y)$ which  transforms \eqref{nfm} to \eqref{nfmt}; thus we have the identity
\[
DH(y)\,\widetilde F(y)=F(H(y)).
\]
Now, with Lemma \ref{lowlem} one has {
\[
DH(y)\widetilde F(y)=\lambda \odot  y+G_*(y)+\widetilde G(y) +Dh(y)\, (\lambda \odot y) + Dh(y)\,G_*(y)+\cdots,
\]
and 
\[
\begin{array}{rcl}
F(y+h(y))&=& F(y)+DF(y)\,h(y)+\cdots \\
              &=& \lambda \odot y+ G_*(y)+F_*(y) + \lambda\odot h(y)+DG_*(y)\,h(y)+\cdots,
\end{array}
\] 
}
with the dots indicating terms of order $\geq 2m$. 
Setting  both sides equal yields 
\[
\left[\lambda,\,h\right]+ {\rm Pr}_{2m-1}\left(\left[G_*,\,h\right]\right)={\rm Pr}_{2m-1}\left(F_*\right)-\widetilde G,
\]
and restricting to the $\delta$-eigenspace yields \eqref{ehds}.
\end{proof}

{
\begin{remark}\label{hdhat}\label{rem:hatnohat} In his convergence proof, Bruno works with the formal power series vector field $\widehat h_\delta$, defined by
    \begin{equation}
\delta\,\widehat  h_\delta+\left[ G_*, \,\widehat  h_\delta\right]=
F_{*,\delta}. 
\end{equation}
In view of \eqref{ehds}, one verifies the alternative representation
\[
h_\delta={\rm Pr}_{2m-1}(\widehat h_\delta).
\]
\end{remark}
}
{
\subsection{Bruno vector fields: A general perspective}
The importance of the representations \eqref{B_or} (equivalently,  \eqref{sys}) seems not to have been fully understood for a long time and almost has not been used by other authors. One  possible reason is  that Bruno worked with the vector space of fields of the form \eqref{B_or} but did not notice a Lie algebra structure on these vector fields.  Only recently it was observed in  \cite{PR} that Bruno vector fields are not merely ordinary vector fields written in an unusual notation, but can be interpreted as a particular class of {discrete vector fields} admitting a Lie algebra structure. Namely, in \eqref{sys} one may regard each point
\[
Q\in{\bf N}\subset\mathbb Z^n
\]
as carrying a vector
$
F_Q\in\mathbb C^n.
$
Thus a Bruno vector field may be viewed as an assignment
\[
Q\longmapsto F_Q
\]
from a discrete subset of the lattice \(\mathbb Z^n\) to \(\mathbb C^n\), together with the reconstruction rule
\[
(F_Q,Q)\longmapsto (x\odot F_Q)x^Q.
\]

This interpretation naturally leads to a more general construction,  introduced in \cite{PR} under the name
{\em{semilattice vector field}},
where  the dimension of the lattice on which the coefficients are indexed can be different from   the dimension of the vectors attached to its points. 
An important feature of this construction is that that it is possibly to introduce  a Lie algebra structure on this set which in such a way that it agrees with usual Lie algebra of vector fields in the case $n=m$.  

}


\section{A convergence theorem}\label{sec:converge}
In this section we consider $F$ is as in \eqref{Asn}, but we will assume that the matrix $A$ is diagonal, 
\begin{equation}\label{eq:Ass} A=\diag(\lambda),\quad\text{with  }\lambda=(\lambda_1,\ldots,\lambda_n)^\top; \quad Ax=\lambda\odot x.
\end{equation}
Moreover we will assume that $F$ is convergent in some neighborhood of $0$.

Generally, convergence criteria for normalizing transformations involve a diophantine condition on the eigenvalues of $A$, and an algebraic condition on the shape of the normal form. We consider these one by one.
\subsection{Diophantine conditions}
\begin{defin} 
For every positive integer $k$ let 
\begin{equation}
    \omega_k:=\min  \left\{ | \la Q,\lambda\ra |;   \ Q\in {\bf N},
\, \la Q,\lambda\ra \ne 0,   \|Q\|< 2^k\right\}. 
\end{equation}
One says that Bruno's {\em Condition } $\omega$ is satisfied if and only if
\begin{equation}
    \sum_{k=1}^\infty \dfrac{\log(\omega_k^{-1})}{2^k}<\infty.
\end{equation}
\end{defin}

\begin{remark}\label{rem21}
An earlier (more restrictive) diophantine condition is due to Siegel \cite{Siegel}, with a later generalization by Pliss \cite{P}. It requires the existence of positive constants $C$ and $\nu$ such that
\[
|\left<Q,\lambda\right>|\geq C\cdot |Q|^{-\nu}\quad\text{  whenever  }   Q\in {\bf N},  \|Q\|\ge 1 ,    \,\left<Q,\lambda\right>\not=0.
\]
Since $|Q|\le \| Q\|+2\le 2^{k+1}$  whenever $\|Q\|< 2^k$, we see that
\[
-\log\omega_k\leq -\log C+\nu (k+1)\log 2,
\]
whenever the Siegel-Pliss condition holds, thus Condition $\omega$ is clearly satisfied.
\end{remark}
Both Condition $\omega$ and the Siegel-Pliss condition may be difficult to check directly for given $\lambda$. But one general result, that is worth mentioning, does not seem to be available in the literature.\footnote{{Siegel may well have been aware of this fact, concerning his original condition.}} 
\begin{teo}
    \label{prop:omega}
Given  $A$ as in \eqref{eq:Ass}, assume that all eigenvalues are algebraic numbers. Then the Siegel-Pliss condition is satisfied.
\end{teo}
\begin{proof}
 For $Q=(q_1,\ldots,q_n)\in\mathbb Z^n$ we abbreviate
\[
\left<Q,\,\lambda\right>:=\sum q_i\lambda_i\text{  and  } |Q|:=\sum|q_i|.
\]
We will show that there exist positive $C$ and $\nu$ such that 
\begin{equation}\label{Qcond}
|\left<Q,\lambda\right>|\geq c\cdot |Q|^{-\nu}\text{  for all  } Q\in \mathbb Z^n\text{  with  } |Q|\geq 2.
\end{equation}
Since for $Q=(p_1,\ldots,p_n)-(\delta_{ij})_{1\leq i\leq n}$ one has $|Q|\leq 1+\sum p_i$,
 the assertion follows (with modified constants).\\
Recall that a complex number is called algebraic when it is a root of a monic polynomial with rational coefficients, and an algebraic integer when it is a root of a monic polynomial with integer coefficients. (We refer to Lang \cite{Lang}, Ch.~VI for pertinent algebraic facts.) 
\begin{enumerate}[(i)]
\item By the hypothesis, there is a finite Galois extension $L$ of $\mathbb Q$ that contains all $\lambda_i$. Let $G$ be the Galois group of $L$, and set $|G|=\nu+1$. We first prove \eqref{Qcond} for the case that all $\lambda_i$ are algebraic integers. Consequently, for any $Q=(q_1,\ldots,q_n)\in \mathbb Z^n$, one has that $\left<Q,\,\lambda\right>$ is an algebraic integer. The norm of an algebraic integer, i.e., the product of all its conjugates under the Galois group, is a rational integer (a power of the constant coefficient of its minimal polynomial, up to sign.) Thus we have
\[
\prod_{\sigma\in G}\left(\sum_{1\leq i\leq n} q_i\sigma(\lambda_i)\right)\in \mathbb Z.
\]
\item 
In case $\left<Q,\,\lambda\right>\not=0$ this implies
\[
1\leq |\left<Q,\,\lambda\right>|\cdot \prod_{\sigma\in G\setminus\{\rm id\}}\left|\sum_i q_i\sigma(\lambda_i)\right|.
\]
With 
\[
C^*:= \max\left\{|\sigma(\lambda_i)|; \, \sigma \in G,\,1\leq i\leq n\right\}
\]
one has
\[
\left|\sum_i q_i\sigma(\lambda_i)\right|\leq C^*\cdot |Q|
\]
for every $\sigma\in G$, and therefore
\[
1\leq |\left<Q,\,\lambda\right>| \cdot {C^*}^\nu\cdot |Q|^\nu
\]
for all $Q\in\bf N$. This is the assertion when all $\lambda_i$ are algebraic integers.
\item Turning to general algebraic numbers, there exists a natural number $N$ such that all $N\cdot\lambda_i$ are algebraic integers. From the reasoning in part (ii) we then obtain 
\[
\dfrac1N \leq |\left<Q,\,\lambda\right>| \cdot {C^*}^\nu\cdot |Q|^\nu
\]
for all $Q\in \bf N$, and the proof is finished.
\end{enumerate}
\end{proof}
\begin{remark}
 Special instances of the theorem can be found in the literature; see e.g.\ 
Proposition 3.2 in the recent paper by Huang et al.\ \cite{HRZ}. 
\end{remark}
One convergence result is a direct consequence of Theorem \ref{prop:omega}:
\begin{cor}
Assume that $Ax=A_sx=\lambda\odot x$, and that the eigenvalues of $A$ are algebraic. If $\la P,\lambda\ra=0$, $P\in M_\lambda$, holds only for $P=0$, then there exists a convergent transformation from  \eqref{sys}  to its linearization $\dot x=Ax$.
\end{cor}
\begin{proof}
{The eigenvalue condition implies that any normal form is necessarily linear and semisimple, hence Bruno's Simplified Condition A, which generally requires that a normal form is of the shape 
\be \label{BrunoSCA} \dot x =(1+s(x))A_sx\ee for some formal power series $s$, $s(0)=0$ (see \cite{Bru}, Ch.~III, \S3) holds, and this implies convergence. }
\end{proof}
\begin{remark}
As stated above, the diophantine criteria by Bruno and Siegel seem to require explicit knowledge of the eigenvalues. But explicit knowledge is not always necessary. 
As a case in point, for any matrix $B\in\mathbb C^{n\times n}$ with algebraic entries, the eigenvalues of $B$ satisfy the Siegel-Pliss condition. 
This holds because the entries of $B$ are contained in some finite extension $K$ of the rational number field $\mathbb Q$ and, in turn, the eigenvalues are contained in some finite extension of $K$, thus are algebraic over $\mathbb Q$. This observation is relevant in the coordinate-independent version of normal forms, where the condition amounts to $[A_sx,F(x)]=0$. (Note that for a given problem one may not be able to explicitly transform a matrix to Jordan canonical form, so coordinate-independent settings are of practical interest.)
\end{remark}

\subsection{Algebraic conditions}\label{subsec:algcond}
To properly state the convergence conditions involving the shape of the normal form, we consider additive decompositions of $\lambda$.
\begin{defin}\label{def:lambdadecomp}
    Let $r\geq 1$ and $\lambda^{(1)},\ldots,\lambda^{(r)}\in \mathbb C^n$ be linearly independent, and define diagonal matrices $A_j$ by $A_j\,x=\lambda^{(j)}\odot x$.
    \begin{enumerate}[(i)]
        \item We call $(\lambda^{(1)},\ldots,\lambda^{(r)})$ an {\em additive decomposition} of $\lambda$ if there exist constants $\gamma_1,\ldots,\gamma_r$ such that 
        \be \label{ladec}
        \lambda=\sum \gamma_j\lambda^{(j)}.
        \ee 
        \item Adopting terminology from Stolovitch \cite{St}, we say that the {set $\{\lambda^{(1)},\ldots,\lambda^{(r)}\}$ } is {\em isoresonant} with respect to $\lambda$ if every polynomial first integral of $Ax$
        is also a first integral of all $A_jx$, $1\leq j\leq r$. Equivalently, $\left<P,\lambda\right>=0$ implies that every $\left<P,\lambda^{(j)}\right>=0$ for all $P\in \mathbb Z_+^n$.
        {We say that a vector field \eqref{sys} in normal form, with an additive decomposition \eqref{ladec}, satisfies the isoresonance condition if $\{\lambda^{(1)},\ldots,\lambda^{(r)}\}$ is isoresonant with respect to $\lambda$.}
        \item Adopting and modifying terminology from Stolovitch \cite{St}, we say that $\lambda^{(1)},\ldots, $ $\lambda^{(r)}$ {\em lie in the diophantine hull} of $\lambda$ if there exists a constant $c>0$ such that
        {
        \begin{equation}\label{eq:DioHull}
           | \left<Q,\lambda^{(j)}\right>|\leq c\cdot |\left<Q,\lambda\right>|, \quad 1\leq j\leq r,\text{  for all  } Q\in{\bf N}.
        \end{equation}
        }
        {We say that a vector field \eqref{sys} in normal form, with an additive decomposition \eqref{ladec}, satisfies the diophantine hull condition if $\lambda^{(1)},\ldots,\lambda^{(r)}$ lie in the diophantine hull of $\lambda$.}
    \end{enumerate}

\end{defin}
\begin{remark}\label{rem26}
    By these definitions, our version of the diophantine hull condition implies the isoresonance condition. 
\end{remark}

\begin{defin}\label{def:14}
{
     Let $G(x)=Ax+\cdots$ be in normal form, with an additive decomposition \eqref{ladec} that satisfies the diophantine hull condition. We say that this system satisfies {\em Condition AS} if
             \begin{equation}\label{GAS}
 G(x)
 =\sum_{P\in M_\lambda^+} (x \odot   G_{P}) x^{P},
        \end{equation}
        with 
        \begin{equation}\label{GASP2}
            G_P=\sum_j \beta_{j,P}\lambda^{(j)};\quad \beta_{j,P}\in \mathbb C,\text{  for all  }P.
        \end{equation}
    }
     \end{defin}
\begin{remark}\label{ASrem} 
    \begin{itemize}
    \item { One may restate condition \eqref{GAS} as 
     \be \label{GAS5}
    G(x)=\sum_{j=1}^r (\gamma_j+\check s_j(x))\cdot \lambda^{(j)}\odot x=\sum_{j=1}^r (\gamma_j+\check s_j(x))\cdot A_jx
    \ee   
    with 
    formal power series $$\check s_j = \sum_{Q\in M_\lambda^+} s_j^{(Q)} x^Q, $$  
    all $\check s_j(0)=0$, $j=1,\dots,r$.
    }
        \item {The normal form property and the diophantine hull condition imply that that every $\check s_j$ is a first integral of $\dot x=A_k x$, for all $k$. Indeed, for $1\leq j,\,k\leq r$ one finds
        \[
        \left[ \lambda^{(k)}\odot x, \,x^Q\,\lambda^{(j)}\odot x\right] =\la Q,\,\lambda^{(k)}\ra\,x^Q\,\lambda^{(j)}\odot x.
        \]
        By the normal form property and a similar computation we have $\la Q,\,\lambda\ra=0$, hence $\la Q,\,\lambda^{(k)}\ra=0$ by \eqref{eq:DioHull}.}
        \item By the characterization of transformations mapping normal forms to normal forms (Bruno \cite{BruB}, \S1, Theorem 2, or also \cite{WalcherNF}, Proposition 1.5), the expansion of such a transformation contains only resonant terms. Using this fact one verifies:  Whenever Condition AS holds for one normal form of a vector field then it will hold for all its normal forms. 
    \end{itemize}
\end{remark}

{The name ``Condition AS'' is chosen because it combines Bruno's Condition A with some notions introduced by Stolovitch.}
{
\begin{exstar}
    \begin{itemize}
        \item With $r=1$ we have an additive decomposition $\lambda=\lambda^{(1)}$, which obviously satisfies the diophantine hull condition. Here Condition AS is Bruno's simplified Condition A as defined in \eqref{BrunoSCA}.
        \item Let $\lambda=\lambda^{(1)}$ and $\lambda^{(2)}=\bar \lambda$ be linearly independent. Then we have an additive decomposition $\lambda=1\cdot\lambda^{(1)}+0\cdot \lambda^{(2)}$, which satisfies the diophantine hull condition due to $\left<Q,\lambda^{(2)}\right>=\overline{\left<Q,\lambda\right>}$. Here Bruno's Condition $A_2$ (see \cite{BruB} \S3) implies Condition AS; see subsection \ref{subsmoremero} below for a further discussion.
    \end{itemize}
\end{exstar}
}

\subsection{A convergence theorem}
The following convergence theorem is an extension of a result of Bruno (with some additions motivated by Stolovitch); our proof 
follows closely Bruno's original reasoning. We use, however, a different formalism that makes the proof simpler and more transparent. 
\begin{teo}\label{BSH}
    Let the differential system \eqref{Asn} be given, with analytic $F$, and assume that
    \begin{itemize}
        \item the eigenvalues of $A$ satisfy Condition $\omega$; and
        \item some normal form of $F$ satisfies Condition $AS$.
    \end{itemize}
    Then there exists a convergent normalizing transformation.
\end{teo}

{For the cases of Condition AS with  $r=1$ and $r=2$ we obtain:
\begin{cor}
\begin{itemize}
\item  If Bruno's Simplified Condition A (condition (\ref{BrunoSCA})) holds, then there exists a convergent normalizing transformation.
\item If Bruno's Condition $A_2$ holds, then there exists a convergent normalizing transformation.
 \end{itemize}
\end{cor}
} 

In our proof we will mimic the proof in Bruno \cite{BruB} for the case $A_2$ (Theorem 4 on p.\ 186). We first sketch a ``blueprint'' for the line of his arguments:
\begin{itemize}
    \item Iteration step, ``blockwise normalization''  { or Newton-Pliss iterations  (see \cite{BruB}, Lemma 3 on p.\ 188 and \cite{P})}: There exists a transformation from a system that is in normal form up to order $m-1$ to a system that is in normal form up to order $2m-1$. This is the content of our Proposition \ref{pro:1n}. This statement holds generally, with no restrictions on $A$ or the shape of the normal form. (Lemma 4 on p. 190 of \cite{BruB} states that Condition A continues to hold for the transformed vector field. Compare Remark \ref{ASrem}.)
    \item Lemma 5 in \cite{BruB}, p.\ 192, yields basic estimates for the transformed vector field and the transformation in each iteration step. Here, the nilpotency of a certain operator plays a crucial role. (Condition AS still ensures nilpotency).
    \item With the above, one sees that the hypotheses of the Corollary in Bruno \cite{BruB}, 
    p.\ 204 are satisfied. Then the remaining steps of the proof are, mutatis mutandis, identical to the steps $1^o$ to $4^o$ in Bruno's proof.
\end{itemize}
\color{black}


\subsection{The homological equation with Condition AS}
We now assume the hypotheses of Theorem \ref{BSH}, and use the notation introduced in the previous section. First we take a closer look at the iteration step, with the help of Lemma \ref{lem:Deltasq}.


\begin{lem}\label{lem:26}
 If the vector field  
 $$
 G(x)
 =\sum_{P\in M_\lambda } (x \odot   G_{P}) x^{P},
 $$
 satisfies condition  \eqref{GAS}, then $\dg^2=0$; in particular $\dg$ is nilpotent.
\end{lem}
\begin{proof}
 If  condition \eqref{GAS} holds, then for any term $(x \odot G_P) x^P$ of $G$ we have that  $G_P$ has the form \eqref{GASP2}. By the isoresonance condition, this yields  $\la P, G_P\ra=0, $ 
  so using  \eqref{DG2} we conclude that $\Delta^2_G=0.$ 
\end{proof}

\begin{lem}\label{lem:hD_n} 
Let $\kappa=(\kappa_1,\dots,\kappa_n)^\top \in \C^n$ 
and 
\[
W=\sum_{Q\in{\bf N}: \la Q,\kappa\ra=\delta} (x\odot W_Q) \,x^Q.
\]
Then:
\begin{enumerate}[1)]
\item $ [W, \kappa\odot x]= -\delta W$.
\item For 
$$
V=(x\odot \kappa) s(x),
$$
where $s(x)=\sum_{P\in M_\kappa} s_P x^P,$ we have
 \be \label{Dhd}
\Delta_{W} V= \delta  s(x)  W.
\ee
\end{enumerate}
\end{lem}
\begin{proof}
1) By Lemma \ref{lem:9} and  \eqref{DelU}  we have
$$
[W,\kappa\odot x]= \Delta_{\kappa} W-\Delta_{W} \kappa\odot x=-\sum_Q
\la Q,\kappa \ra  (x\odot W_Q) x^Q=-\delta \sum_Q (x\odot W_Q) x^Q =-\delta   W. 
$$

2) By \eqref{DelU} 
$$
\Delta_{W} V=\sum_{P,Q} \la Q, \kappa\ra s_P (x\odot W_Q) x^{P+Q}= \delta  s(x) W. 
$$ 
Compare also the third item in Lemma \ref{linaclem}.
\end{proof}

{

\begin{lem}\label{lem_solhomol_new}
Let  $\lambda^{(1)}, \dots, \lambda^{(r)}\in \C^n$ be linearly independent and $\lambda=\gamma_1\lambda^{(1)}+\cdots+\gamma_r\delta^{(r)}$ an additive decomposition of $\lambda$. 
For $Q\in\bf N$ and $\delta=\left<Q,\,\lambda\right>$ write
\[
\delta^{(i)}=\left<Q,\,\lambda^{(i)}\right>; \quad 1\leq i \leq r.
\]
{Assume that $G(x)$ satisfies Condition AS with \eqref{GAS}. Then, with the notation of Proposition \ref{pro:1n},
\[
G_*(x)=\sum_{j=1}^r (\gamma_j+s_j(x))\cdot \lambda^{(j)}\odot x=\sum_{j=1}^r (\gamma_j+ s_j(x))\cdot A_jx
\]
    with suitable polynomials\footnote{ Note the contrast to Definition \ref{def:14}, where formal  power series are involved.  }
     $ s_1,\ldots, s_r$, such that all $ s_j(0)=0$, and
 the unique solution to \eqref{ehds} is 
\begin{equation}\label{eq:58}
         h_{\delta}(x)= {\rm Pr}_{2m-1} \left(  \dfrac{1}{ (\delta+ \sum \delta^{(j)}s_j(x))}\left({\rm Id}+\dfrac{1}{ (\delta+ \sum \delta^{(j)} s_j(x))}\Delta_{G} \right) F_{*, \delta}(x)\right).
    \end{equation}
}
    where  $\delta=\delta^{(1)} \gamma_1 +\cdots+\delta^{(r)} \gamma_r $.
\end{lem}
\begin{proof}
Since for any $U,V,W\in \mc M$ it holds that
$$
\Delta_{U+V}W=\Delta_UW+\Delta_VW,
$$
we can write 
$$
\Delta_{G_*}=\sum_{i=1}^r\Delta_{G_{*,i}}; 
$$
where 
\be \label{Gil}
G_{*,i}(x)= (x \odot \lambda^{(i)}) (\gamma_i+ s_i(x)),
\ee
and  $G_*=\sum_{i=1}^r G_{*,i}$ by \eqref{GAS}.

{We look for the solution 
$$\widehat h_\delta= \sum_{\substack{Q:\,\la Q,\lambda\ra=\delta\\
                \|Q\|\ge m}}(x\odot\widehat  H_Q) x^Q
$$ 
of equation \eqref{hdhat}, noting from Remark \ref{rem:hatnohat} that $h_\delta={\rm Pr}_{2m-1}(\widehat h_\delta)$
solves \eqref{ehds}.}

In view of Lemma \ref{lem:hD_n}
the computations yield 
$$
\begin{aligned}
[\widehat h_{\delta}, G_*]=& \sum_{i=1}^r [\widehat h_{ \delta}, (x\odot  \lambda^{(i)})(\gamma_i+ s_i(x))] =\\
&\sum_{i=1}^r \left( -\delta^{(i)}\gamma_i  \widehat h_{\delta} -\Delta_ {\widehat h_{\delta}} \widetilde G_{*,i}+\Delta_{\widetilde G_{*,i}}\widehat  h_\delta\right) =
\sum_{i=1}^r\left(  -(\delta^{(i)}\gamma_i +\delta^{(i)}  s_i(x)) {\rm Id} +\Delta_{\widetilde G_{*,i}}\right)\widehat  h_{ \delta},
\end{aligned}
$$
where $\widetilde G_{*,i}= G_{*,i}- \gamma_i (x\odot \lambda^{(i)})$.
Therefore, taking into account that $\Delta_{G_{*,i}}=\Delta_{\widetilde G_{*,i}}$, we can write \eqref{hdhat} as 
 \be \label{ASeq}
   \left(\sum_{i=1}^r (\gamma_i\delta^{(i)}+\delta^{(i)}s_i(x)) {\rm Id}-\Delta_{G_{*,i}}\right)\widehat  h_{ \delta} = F_{ *,\delta}.  
\ee
  Inverting the geometric series, using  the nilpotency of $\Delta_{G_{*,i}}$ (by Lemma \ref{lem:26}), and  the fact that $\Delta_{G_*}=\sum_{i=1}^r\Delta_{G_{*,i}}$
we obtain
$$
  \widehat h_{\delta}(x)=    \dfrac{1}{ (\delta+ \sum \delta^{(j)}s_j(x))}\left({\rm Id}+\dfrac{1}{ (\delta+ \sum \delta^{(j)} s_j(x))}\Delta_{G} \right)  F_{*, \delta}(x),
  $$
  and we see that \eqref{eq:58} is a solution to \eqref{ehds}.
\end{proof}

{
Setting $\alpha^{(i)}=\delta^{(i)}/\delta, \ (i=1,\dots, r)$ we can rewrite \eqref{eq:58}
as 
 \begin{equation}\label{eq:invert_delta_b}
        h_{\delta}(x)=   {\rm Pr}_{2m-1} \left( \dfrac{1}{ \delta (1 +\sum \alpha^{(j)} s_j(x))}\left({\rm Id}+\dfrac{1}{ \delta(1+\sum \alpha^{(j)}s_j(x))}\Delta_{G_*}\right) F_{*,\delta}(x)\right).
    \end{equation}
}
 }

\begin{remark}
    Up to this point, the proofs require only a normal form \eqref{GAS} such that the system $(\lambda^{(1)},\ldots,\lambda^{(r)})$ is isoresonant with respect to $\lambda$.
\end{remark}


\subsection{Estimates}\label{A2est}
We recall some notation; compare Bruno \cite{BruB,Bru}.
\begin{defin}\label{def:normstuff}
    \begin{enumerate}
        \item For a series 
        \[
        a(y)=\sum a_Q y^Q,\quad a_Q\in\mathbb C,
        \]
        and a real number $\rho>0$ we set
        \[
        \overline{| a(y)|}:= \sum |a_Q|y^Q, \quad \overline{ |a|}_\rho:= \sum |a_Q|\rho^{\|Q\|}.
        \]
        \item Given a series $b(y)=\sum b_Qy^Q$ ($b_Q\ge 0$), we write
        \[
        a(y)\prec b(y) \Longleftrightarrow |a_Q|\leq b_Q,\quad\text{for all  }Q.
        \]
        \item For $B=\begin{pmatrix} b_1\\ \vdots \\b_n\end{pmatrix}$ we set
        \[
        |B|:= |b_1|+\cdots+|b_n|; \quad |B|_\rho:= |b_1|_\rho+\cdots+|b_n|_\rho; 
        \]
        and similarly for matrices.
        \item For $B=\begin{pmatrix} b_1\\ \vdots \\b_n\end{pmatrix}$ we set
        $
        \overline{|B|}:=\begin{pmatrix} \overline{|b_1|}\\ \vdots \\ \overline{|b_n|}\end{pmatrix}; 
        $
         and similarly for matrices.

    \end{enumerate}
\end{defin}

\begin{lem}\label{lem:stepest}
Consider the setting of Proposition \ref{pro:1n} and Corollary \ref{cor47}.
Assume that for some $\rho$, $\frac12<\rho\leq 1$, one has 
\be \label{est}
 \overline { |F|}_\rho<1, \ |\overline{ G_*}|_\rho <c_1,   \ |\overline{ \Delta_{G_*}}|_\rho <c_1,
\ee
with some constant $c_1$, and moreover assume that Condition AS holds with linearly independent $\lambda^{(1)}, \dots, \lambda^{(r)} $.  Furthermore let $m=2^k$. For $Q\in{\bf N}$ with $m\leq \Vert Q\Vert\leq 2m-1$ let $\delta=\left<Q,\,\lambda\right>$, $\delta^{(i)}=\left<Q,\,\lambda^{(i)}\right>$. 
Then there exists a constant $c_2$ such that 
\be \label{est_h_bn}
\overline{ |h|}_\rho <\dfrac{c_2}{\omega_{k+1}^2}.
\ee
\end{lem}
\begin{proof}
We abbreviate $\alpha^{(i)}=\delta^{(i)}/\delta$, as in \eqref{eq:invert_delta_b}. Note also that since $\lambda^{(1)}, \dots, \lambda^{(r)} $ satisfy Condition AS
 there is a constant $c$
 depending only on $\lambda$ such that,  for all $Q$,
\be \label{eq:c_est}
|\delta^{(i)}/\delta|<c.
\ee 
\begin{enumerate}[(i)]
\item
We first show that there is a $c_1$ such that under the conditions of the lemma
\be \label{eq:asr}
|\alpha^{(i)}|  \overline{|s_i|}_\rho < \frac{1}{2r }.
\ee
We write each $s_j(x)$ as the series 
$$
s_j(x)=\sum_{P} s^{(j)}_P x^P.
$$
By \eqref{GASP2}
we have for each resonant $P$ 
$$
G_{*,P}=\sum_{i=1}^r \lambda^{(i)} s_P^{(i)}. 
$$
Now the matrix $\Lambda$ with columns $\lambda^{(1)},\ldots,\lambda^{(r)}$ defines a linear map from $\mathbb R^r$ to $\mathbb R^n$, which is injective by the linear independence of the $\lambda^{(i)}$. Denote by $\Vert \cdot \Vert_1$ the 1-norms on the respective spaces.  By compactness, $\Vert\Lambda v\Vert_1$ attains a positive minimum $\beta$ on the unit sphere in $\mathbb R^r$, and therefore $\Vert\Lambda v\Vert_1\geq \beta\,\Vert v\Vert_1$ for all $v\in \mathbb R^r$. We thus find, for $1\leq j\leq r$,
\[
|\alpha^{(j)}|\,|s_P^{(j)}|\leq c\sum_i |s_P^{(i)}|\leq \frac{c}{\beta}\,\Vert G_{*,P}\Vert_1,
\]
using \eqref{eq:c_est}. Altogether we have 
$$
|\alpha^{(i)}| |s_i| \prec \frac{c} \beta |\overline{G}|.
$$
We take 
$$
c_1=\beta\frac{1}{2 r c },
$$
then \eqref{eq:asr} holds and 
$$
|\sum_{j=1}^r \alpha^{(j)}  s_j(x)|_\rho\le
\sum_{j=1}^r |\alpha^{(j)}| | s_j(x)|_\rho\le \frac 12.\
$$
\item
Thus  we have
$$
\left|\frac{1}{ 1+\sum\alpha^{(j)}  s_j(x)}\right|_\rho \le \left|\frac{1}{ 1-|\sum\alpha^{(j)}  s_j(x)|}\right|_\rho  
\le \left|\frac{1}{ 1-\sum |\alpha^{(j)} s_j(x)|}\right|_\rho \le \frac{1}{1-\frac 12} =2.
$$

{
  From \eqref{eq:invert_delta_b} one obtains 
  \[
|\overline{ \widehat h_{ \delta}|}_\rho\leq  \dfrac{2}{|\delta|}\left(n+c_1\dfrac{2}{|\delta|}\right) |\overline{F_{*, \delta}|}_\rho.
\]
Since for terms  of $h_\delta$
 $m\leq \Vert Q\Vert\leq 2m-1$
by the definition of  $\omega_{k+1}$ we have $|\delta| \ge  \omega_{k+1}$. }
By Corollary \ref{cor47}
$
h=\sum_{\delta\ne 0} h_\delta
$,
 where $h_\delta$ is the unique solution to \eqref{ehds} in $\mathcal V^{2m-1}$.  
{Therefore, from \eqref{eq:invert_delta_b} for $h_\delta$ we obtain the estimate
\[
 \overline{ | h_\delta|}_\rho \leq  \overline{ | \widehat h_\delta|}_\rho \leq\dfrac{2}{\omega_{k+1}}\left(n+c_1\dfrac{2}{\omega_{k+1}}\right) |\overline{F_{*, \delta}|}_\rho.
\]
}
Now recall $\omega_{k-1}\leq \omega_1$, and let 
\[
c_2\geq 2\left(n\omega_1+2c_1\right).
\]
Then 
\[
|\overline{  h_{ \delta}|}_\rho\leq \dfrac{ c_2}{\omega^2_{k+1}} |\overline{F_{*, \delta}|}_\rho.
\]
Summing up over all $\delta$ and taking into account that  $|\overline { F_*}|_\rho|\leq \overline { F}|_\rho<1$ we see that  
the above estimate implies \eqref{est_h_bn}.
\end{enumerate}
\end{proof}
\begin{remark}\label{remDvsDelta}
The conditions in equation \eqref{est} are different from those given in Bruno \cite{BruB}, which require
    \be \label{estD}
  |\bar G_*|_\rho <c_1,   \ |\overline{ D{G_*}}|_\rho <c_1.
\ee
But, due to \eqref{DvsDelt}, one sees that
\[
|\overline{ D{G_*}}|_\rho \leq  |\bar G_*|_\rho+ |\overline{ \Delta_{G_*}}|_\rho,
\]
thus \eqref{est} implies $|\overline{ D{G_*}}|_\rho  <2 c_1$, so \eqref{estD} holds with a constant $\widetilde c_1:=2c_1$. The converse estimate is shown from $|\overline{ \Delta_{G_*}}|_\rho\leq |\bar G_*|_\rho+ |\overline{ D{G_*}}|_\rho$. Thus the conditions are indeed equivalent.
\end{remark}
{
\begin{proof}[Proof of Theorem \ref{BSH}]
    The remainder of the proof consists of appealing to Bruno \cite{BruB}: We have shown that the statement of Bruno's Corollary on p.\ 204 of \cite{BruB} holds. From here on, the proof is identical (modulo notation) to the proof of \cite{BruB}, Theorem 4 on p.\ 204 ff.
\end{proof}
}

\subsection{On the diophantine hull condition}

Condition \eqref{eq:c_est} is crucial in the line of proof. Adopting terminology in Stolovitch \cite{St}, we say that the $\lambda^{(i)}$ are in the {\em diophantine hull} of $\lambda$ if \eqref{eq:c_est} is satisfied. We now take a closer look at the diophantine hull condition, which turns out to be quite restrictive. First, we note that it is not satisfied even in seemingly simple situations.
\begin{exstar} {\em 
    Let $\mu>0$ be an irrational number, $\lambda^{(1)}=(1,-1,0,0)$, $\lambda^{(2)}=(0,0,\mu,-\mu)$ and $\lambda=\lambda^{(1)}+\lambda^{(2)}$. Then for any $Q\in{\bf N}$ we have 
    \[
|\left<\lambda^{(1)},Q\right>|\geq 1 \text{ when }\left<\lambda^{(1)},Q\right>\not=0;\quad |\left<\lambda^{(2)},Q\right>|\geq \mu \text{ when }\left<\lambda^{(2)},Q\right>\not=0.
    \]
    But $|\left<\lambda,Q\right>|= |q_1-q_2+\mu(q_3-q_4)|$ may be nonzero and arbitrarily small, as known from approximation of irrational numbers by fractions; see e.g. Hurwitz \cite{Hur}. Therefore the diophantine hull condition does not hold.}
\end{exstar}
To obtain some positive results, we identify $\mathbb C$ with $\mathbb R^2$, equipped with the standard scalar product. 
 \begin{pro}\label{pro:cones} 
     Let $K_1$, $K_2$ be different straight lines through $0$ in $\mathbb C$, and furthermore $\lambda=\lambda^{(1)}+\lambda^{(2)}$ such that every entry of $\lambda^{(j)}$ is contained in $K_j$, $j=1,2$. Then there exists a constant $c>0$ such that, for any $Q\in {\bf N}$ with $\Vert Q\Vert=m\geq 2$ and $\left<\lambda,Q\right>\not=0$ one has 
     \[
     |\left<\lambda^{(1)},Q\right>|\leq c\cdot|\left<\lambda,Q\right>| \text{  and  }|\left<\lambda^{(2)},Q\right>|\leq c\cdot  |\left<\lambda,Q\right>|.
     \]
     Therefore, if $\lambda$ satisfies Condition $\omega$, and a normal form of the analytic vector field $F(x)=\lambda\odot x+\cdots$ satisfies Condition AS, then there exists a convergent transformation of $F$ to normal form.
 \end{pro}
 \begin{proof}
 Here we denote the scalar product in $\mathbb R^2$ by $(\cdot,\cdot)$. By elementary geometry there exists $\rho<1$ such that for any $u_1\in K_1$ and $u_2\in K_2$
 \begin{equation}\label{eq:scapro} 
     |(u_1,u_2)|\leq \rho\,\Vert u_1\Vert \,\Vert u_2\Vert.
 \end{equation}
     Now 
     \[
     \begin{array}{rcl}
       (u_1+u_2,\,u_1+u_2)&\geq    & \Vert u_1\Vert^2+\Vert u_2\Vert ^2-2|(u_1,u_2)| \\
          & \geq &\Vert u_1\Vert^2+\Vert u_2\Vert ^2-2\rho\, \Vert u_1\Vert\,\Vert u_2\Vert \\
         & \geq& (1-\rho) \left(\Vert u_1\Vert^2+\Vert u_2\Vert ^2\right) +\rho\,\left(\Vert u_1\Vert^2+\Vert u_2\Vert ^2-2 \Vert u_1\Vert\,\Vert u_2\Vert\right)\\
         &\geq &(1-\rho) \left(\Vert u_1\Vert^2+\Vert u_2\Vert ^2\right).
     \end{array}
     \]
     Thus 
     \[
     \dfrac{\Vert u_j\Vert ^2}{\Vert u_1+u_2\Vert^2}\leq \dfrac{1}{1-\rho}, \quad j=1,\,2.
     \]
     Applying this inequality to $u_j=\left<Q,\lambda^{(j)}\right>$, the assertion is immediate from Lemma \ref{lem:stepest} and the ensuing arguments.
 \end{proof}
 \begin{remark}{
     In essence we have recovered Bruno's Condition $A_2$: Let $\lambda=\lambda^{(1)}+\lambda^{(2)}$ as above.
     \begin{itemize}
         \item  Let $K_1=\mathbb R \zeta_1$ and $K_2=\mathbb R\zeta_2$ be two different lines such that every entry of $\lambda^{(j)}$ is contained in $K_j$, $j=1,2$.  Since
     \[
     \bar\zeta_2\lambda-\zeta_2\bar\lambda= \left(\zeta_1\bar\zeta_2-\bar\zeta_1\zeta_2\right)\lambda^{(1)}\text{  and  } \bar\zeta_1\lambda-\zeta_1\bar\lambda= \left(\bar\zeta_1\zeta_2-\zeta_1\bar\zeta_2\right)\lambda^{(2)},
     \]
there are complex constants $\theta_{ij}$ so that
\[
\begin{array}{rcl}
\lambda^{(1)}&=&\theta_{11}\lambda+\theta_{12}\bar\lambda\\
\lambda^{(2)}&=&\theta_{21}\lambda+\theta_{22}\bar\lambda.
\end{array}
\]
Therefore any normal form 
     \[
     (1+s_1(x))(\lambda^{(1)}\odot x) +  (1+s_2(x))(\lambda^{(2)}\odot x)= \lambda\odot x  +  s_1(x)(\lambda^{(1)}\odot x) +  s_2(x)(\lambda^{(2)}\odot x) 
     \]
     can also be represented as 
     \[
   \lambda\odot x+  \hat a(x)(\lambda\odot x) +\hat b(x)(\bar\lambda\odot x)
     \]
     with series $\hat a$ and $\hat b$, $\hat a(0)=\hat b(0)=0$, and Bruno's Condition $A2$ is satisfied.
     \item Conversely, if $\lambda$ and $\bar \lambda$ are linearly independent, then setting $\lambda^{(1)}={\rm Re}\,\lambda$ and  $\lambda^{(2)}={\rm Im}\,\lambda$ yields the above scenario with $K_1=\mathbb R$ and $K_2=i\,\mathbb R$.
     \end{itemize}
     }
 \end{remark}{
 Comparing the above to Remark 6.3.6 in Stolovitch \cite{St} (specialized to one dimensional Lie algebras), one may suspect that Condition \eqref{eq:c_est} -- which is crucial for the convergence proof mimicking Bruno -- is valid only in scenarios equivalent to Bruno's $A_2$. We present one more example to show how a deviation from the scenario with two lines may destroy the condition. But by a further example we show that there exist instances when Condition AS holds while $A_2$ does not. }

 \begin{exstar}{\em Let $\alpha, \beta\in \mathbb R$ be irrational, and 
 \[
 \lambda^{(1)}:=\diag(1,-1,0,0,0),\, \ \lambda^{(2)}:=\diag(0,0,i,-i,\alpha+i\beta); \ \lambda:=\lambda^{(1)}+\lambda^{(2)},
 \]
noting that the isoresonance condition (see Definition \ref{def:lambdadecomp}) forces at least one of $\alpha,\, \beta$ to be irrational. Then  
 \[
 \left<Q,\lambda^{(1)}\right>=q_1-q_2,\quad \left<Q,\lambda^{(2)}\right>=q_5\alpha +i(q_3-q_4+q_5\beta),
 \]
 and
 \[
 \left<Q,\lambda\right>= (q_1-q_2+q_5\alpha)+i\,(q_3-q_4+q_5\beta).
 \]
 Here $|\left<Q,\lambda^{(1)}\right>|= |q_1-q_2|$ can be made arbitrarily large, while for given $q_1,\,q_2$ the absolute value of the real part of $\left<Q,\lambda\right>$ can be made $<1$ by suitable choice of $q_5$ (Hurwitz), and for given $q_5$ the absolute value of the imaginary part of $\left<Q,\lambda\right>$ can be made $<1$ by suitable choice of $q_3,\,q_4$.Therefore, the diophantine hull condition cannot hold here.\\
This example extends to dimension $>5$, with appending more zeros to $\lambda^{(1)}$, and arbitrary entries to $\lambda^{(2)}$.}
  \end{exstar}
  In the following example, $A_2$ is not applicable because $\lambda$ and $\bar\lambda$ are linearly dependent.
  \begin{exstar}{\em { 
      Let $\beta_2$ and $\beta_3$ be positive real numbers, moreover
      \[
      \lambda^{(1)}=\begin{pmatrix}
          0\\ \beta_2\\0
      \end{pmatrix},\quad \lambda^{(2)}=\begin{pmatrix}
          0\\ 0\\ \beta_3
      \end{pmatrix}, \quad\text{and}\quad \lambda=\lambda^{(1)}+\lambda^{(2)}.
      \]
      This additive decomposition satisfies the diophantine hull condition: We consider 
      \[
      \la Q,\,\lambda\ra= q_2\beta_2+q_3\beta_3,\quad \la Q,\,\lambda^{(1)}\ra= q_2\beta_2
      \]
      for $Q\in \bf N$, and show that there exists a constant $c_1>0$ so that $ \la Q,\,\lambda\ra\geq c_1\cdot \la Q,\,\lambda^{(1)}\ra$. We may assume that $q_2\not=0$.\begin{itemize}
          \item When $q_2>0$ and $q_3\geq 0$ then $|q_2\beta_2+q_3\beta_3|\geq |q_2\beta_2|$ is obvious.
          \item In case $q_2=-1$ note \[(-\beta_2+q_3\beta_3)/(-\beta_2)= 1-q_3\beta_3/\beta_2.\] There exists $q_3^*$ so that $1-q_3\beta_3/\beta_2>0$ if and only if $q_3<q_3^*$, and the minimum, $\mu_1>0$, is attained at $q_3^*-1$. For $q_3\geq q_3^*$ one has $1-q_3\beta_3/\beta_2\leq 1-q_3^*\beta_3/\beta_2$ and $|1-q_3\beta_3/\beta_2|\geq |1-q_3^*\beta_3/\beta_2|=:\mu_2>0$.
          \item In case $q_3=-1$ note
          \[
          (q_2\beta_2-\beta_3)/(q_2\beta_2)= 1-\beta_3/(q_2\beta_2).
          \]
          There exists $q_2^*$ so that $1-\beta_3/(q_2\beta_2)<0$ if and only if $q_2<q_2^*$. For these finitely many $q_2$ one has a maximum, $-\mu_3$, with $\mu_3>0$. For all $q_2\geq q_2^*$ one finds $1-\beta_3/(q_2\beta_2)\geq1-\beta_3/(q_2^*\beta_2)=:\mu_4>0.  $
      \end{itemize}
      Altogether one sees $\la Q,\lambda\ra\geq c_1\la Q,\lambda^{(1)}\ra$, with the minimum $c_1$ of $1,\,\mu_1,\ldots,\mu_4$. The same argument shows the desired inequality for $\lambda^{(2)}$.
      
      A normal form for this equation (with the diophantine hull condition) reads
      \[
      \dot x=\begin{pmatrix} 0\\ \beta_2\\ \beta_3\end{pmatrix}+ \begin{pmatrix}
          0\\  (1+\check s_2(x_1))x_2\\ (1+\check s_3(x_1)) x_3
      \end{pmatrix}
      \]
      with power series $\check s_i$, $\check s_i(0)=0$.
      
      For this system, convergence also would follow from Theorem \ref{Thm:Bruno1*} in section \ref{sec:Bruno1} below. But here we have a direct (and less involved) proof.}}
  \end{exstar}

\section{Normalization of commuting vector fields}\label{sec:stolo}

In this section we consider the simultaneous normalization of pairwise commuting analytic or formal vector fields $F^{(1)},\ldots,F^{(s)}$ on $\mathbb C^n$, with
\begin{equation}\label{eq:abellie}
F^{(k)}(x)=A^{(k)}x+\text{t.h.o} =\lambda^{(k)}\odot x +\text{t.h.o},
\end{equation}
and linearly independent $\lambda^{(1)},\ldots,\lambda^{(s)}$. One says that this abelian Lie algebra is in normal form if
\begin{equation}
\left[A^{(k)}x,\,F^{(\ell)}(x)\right]=0,\quad\text{ for all  }k,\,\ell.
\end{equation}

The results and the original proofs regarding simultaneous normalization go back to Stolovitch \cite{St}. As will be shown in this section, one can prove results of this type with Bruno's approach, preceded by  a few technical steps. Again we will not discuss the most general setting.

As a preliminary step  we consider transformations (as in subsection \ref{subsec:trafos}) and note that they respect Lie brackets. For the readers' convenience we include a direct  proof of this familiar fact.
\begin{lem}
Let formal vector fields $F_1,\,F_2$ be given, and let $H$ be an invertible formal transformation such that the substitution $x=H(y)$ transforms $F_1$ to $\widetilde F_1$ and $F_2$ to $\widetilde F_2$. Then this substitution transforms $\left[ F_1,\,F_2\right]$ to $\left[\widetilde  F_1,\,\widetilde F_2\right]$.
\end{lem}
\begin{proof}
According to \eqref{sopremap} we have $DH(y)\,\widetilde F_j(y)=F_j(H(y))$ for $j=1,\,2$.  Differentiation yields
\[
D^2H(y)\left(\widetilde F_j(y),\,*\right)+DH(y)\,D\widetilde F_j(y)= DF_j(H(y))\,DH(y),
\]
and furthermore
\[
\begin{array}{rcl}
D^2H(y)\left(\widetilde F_1(y),\,\widetilde F_2(y)\right)+DH(y)\,D\widetilde F_1(y)\widetilde F_2(y)&=&DF_1(H(y))\,DH(y),\widetilde F_2(y)\\
&=&DF_1(H(y))\,F_2(H(y)).\\
\end{array}
\]
Exchanging indices, taking the difference and using the symmetry of the second derivative yields
\[
DH(y)\,[\widetilde F_1,\,\widetilde F_2](y)=[F_1,\,F_2](H(y)),
\]
as asserted.
\end{proof}

In the following we will give a different proof for some results by Stolovitch \cite{St} on formal and convergent simultaneous normalization of commuting vector fields, along the lines of Bruno's arguments. Just as we discussed Bruno's proof only in simplified settings  (Simplified Condition A and Condition AS), we will discuss Stolovitch's proof only in a simplified setting (a variant of Simplified Condition A).
Our  main purpose is to provide a transparent proof that requires no advanced machinery.

\subsection{Formal normalization}

We start by introducing notation, and making some preliminary observations. Let $m\geq 1$ and abbreviate $\mathcal V=\mathcal V^{2m-1}$, the space of polynomial vector fields of order $\leq 2m-1$. Then the linear maps ${\rm ad}_{A^{(k)}}=\left[A^{(k)},\cdot\right]$ are simultaneously diagonalizable on $\mathcal V$, thus
\[
\mathcal V=\bigoplus \mathcal V_{\bar\delta}, \quad \bar\delta=(\delta_1,\ldots,\delta_s),
\] 
with
\[
R\in \mathcal V_{\bar\delta}\Longleftrightarrow\left[A^{(k)}, \,R\right]=\delta_k\,R,\quad 1\leq k\leq s.
\]
We will write $R\in\mathcal V$ as sum of its eigenspace components,
\[
R=\sum_{\bar\delta} R_{\bar\delta}.
\]
Furthermore we note that
\begin{equation}\label{jacmulti}
\left[\mathcal V_{\bar\delta},\,\mathcal V_{\bar \theta}\right]\subseteq \mathcal V_{\bar\delta+\bar\theta},
\end{equation}
a consequence of the Jacobi identity.\\
Our first goal is to state and prove an extension of Proposition \ref{pro:1n}. At the outset we assume that for some $m\geq 1$ we have pairwise commuting vector fields
\begin{equation}\label{liealgvf}
F^{(k)}(x)=A^{(k)}x+G_*^{(k)}(x)+F_*^{(k)}(x)
\end{equation}
with $1\leq k\leq s$, $1\leq \ord G_*^{(k)}\leq m-1$, $\ord F_*\geq m$, which satisfy the partial normal form conditions
\begin{equation}
\left[A^{(k)},\,G_*^{(\ell)}\right]=0,  \quad \text{ all  }k,\,\ell.
\end{equation}
 Moreover we set 
     \begin{equation}\label{equ:Glie}
         G^{(k)}(y)=A^{(k)}y+G_*^{(k)}(y)+\widetilde G^{(k)}(y)+\cdots,
     \end{equation}
     where each  $\widetilde G^{(k)}$ is the resonant component (i.e., the component in $\mathcal V_{\bar 0}$) of ${\rm Pr}\,(F_*^{(k)})$, and the dots indicate terms of order $\geq 2m$.
     We furthermore set
\begin{equation}
{\rm Pr}_{2m-1}\left(F_*^{(k)}\right)= \widetilde G_*^{(k)}+  \widehat F_*^{(k)},\quad \text{  with  } \widehat F_*^{(k)}={\rm Pr}_{2m-1}\left(\sum_{\bar\delta\not=0}F_{*,\bar\delta}^{(k)}\right).
\end{equation}

\begin{pro}\label{pro:lien3}
Given the commuting vector fields \eqref{liealgvf} satisfying the partial normal form conditions,
there exists $h\in \mathcal{M}$  with all terms in $h$ of order between $m$  and $2m-1$, such that 
  the substitution 
  $$
  x=y+h(y)
  $$
simultaneously  transforms them to the vector fields \eqref{equ:Glie}.
\end{pro}

\begin{proof} We call $Q\in\bf N$ resonant if $\left<Q,\,\lambda^{(k)}\right>=0$ for all $k$, and non-resonant otherwise.
We order  all non-resonant $Q$,  $ m\le \|Q\| \le 2m -1$ using the degree lexicographic ordering
and  write   
\be \label{hprmn32}
 \widehat F_*^{(k)}(x)=\sum_{Q} (x\odot  \widehat  F^{(k)}_{*,Q}) x^{Q}
\ee 
(so some $\widehat F^{(k)}_{*,Q}$ may be zero).

Let $Q_1$ be the smallest among the $Q$ with $\widehat F_{*,Q}^{(i)}\ne 0 $ for some $i\in\{1,\dots,s\}$, and call $\delta_k=\left<Q_1,\lambda^{(k)}\right>$ for all $k$. In view of \eqref{Mpq},
evaluating  $\left[F^{(k)},\,F^{(\ell)}\right]=0$ 
and comparing terms of degree $Q_1$ yields
\begin{equation}\label{AFcomp2}
\left[A^{(k)},\widehat F_{*,Q_1}^{(\ell)}\right]=\left[A^{(\ell)},\widehat F_{*,Q_1}^{(k)}\right]
\end{equation}
for all $k,\ell\leq s$, therefore
\begin{equation}\label{AFcompev1}
\delta_k \widehat F_{*,Q_1}^{(\ell)}=\delta_\ell \widehat F_{*,Q_1}^{(k)}.
\end{equation}
Hence, whenever $\delta_k\not=0$ and $\delta_\ell\not=0$, then 
\begin{equation}\label{AFcompev2}
\dfrac1{\delta_\ell}\widehat  F_{*,Q_1}^{(\ell)} = \dfrac1{\delta_k} \widehat F_{*,Q_1}^{(k)}.
\end{equation}

Now the proof closely follows the pattern of the proof of Proposition \ref{pro:1n}.
Fix some $k$ such that $\delta_k\not=0$.
By \eqref{eq:hs2} the transformation 
$$
x=H^{[1]}(y)=y+h_1(y) :=y+\dfrac{1}{\delta_k}(y\odot F^{(k)}_{*,Q_1} )y^{Q_1} 
$$
eliminates the term $ \widehat F_{*,Q_1}^{(k)} $ from the vector field $F^{(k)}$. But actually, due to \eqref{AFcompev2}, this transformation also eliminates the term $ \widehat F_{*,Q_1}^{(\ell)} $ from the vector field $F^{(\ell)}$ whenever $\delta_\ell\not=0$, while $ \widehat F_{*,Q_1}^{(\ell)} =0$ in case $\delta_\ell=0$ by \eqref{AFcompev1}. To summarize, $H^{[1]}$ transforms all $F^{(k)}$ to  commuting vector fields whose non-resonant terms have degree $>Q_1$. 

From here on, the remaining proof steps are the same as in Proposition \ref{pro:1n}, up to obvious modifications of notation (replacing $F$ by $F^{(k)}$, et cetera).
\end{proof}
Furthermore, the proof of Corollary \ref{cor47} directly carries over to the proof of:
\begin{cor}\label{cor4711}
The function  $h$  in Proposition \ref{pro:lien3} can be written as 
$$
h=\sum_{\bar\delta\ne 0} h_{\bar\delta}
$$
where for each $k$ with $\delta_k\not=0$, $h_{\bar\delta}$ is the unique solution of
\begin{equation}\label{ehds2}
\delta_k\,h_{\bar\delta}+{\rm Pr}_{2m-1}\left(\left[ G^{(k)}_*, \,h_{\bar\delta}\right]\right)={\rm Pr}_{2m-1}(F^{(k)}_{*,\bar\delta}) \text{  in  } \mathcal V^{2m-1}.
\end{equation}
\end{cor}
For the sake of completeness we note the result about formal normalization.
\begin{cor}\label{stoloform}
Let $F^{(1)},\ldots,F^{(s)}$ be as in \eqref{eq:abellie}. Then there exists a formal transformation $H$ such that the transformed vector fields $G^{(1)},\ldots,G^{(s)}$ (commute and) are in simultaneous normal form.
\end{cor}
\subsection{A convergence theorem {due to Stolovitch}}
Now we turn to convergence problems, thus we assume that the $F^{(k)}$ are analytic. We state two conditions that go back to Stolovitch \cite{St} and may be seen as extensions of Bruno's Simplified Condition $A$ and Condition $\omega$ {to the Lie algebra setting.}
\begin{defin}
Let $F^{(1)},\ldots,F^{(s)}$ be as in \eqref{eq:abellie}. 
\begin{itemize}
\item We say that some (formal) normal form $G^{(1)},\ldots,G^{(s)}$ satisfies {\em Condition AL}
 if there are scalar formal power series $\check v_{ik}$ with $\check v_{ik}(0)=0$ such that
\begin{equation}\label{eq:ALone}
G^{(k)}(x)=\sum_{i=1}^s\left(\delta_{ik}+\check v_{ik}(x)\right) A^{(i)}x,\quad 1\leq k\leq s;
\end{equation}
here $\delta_{ik}$ denotes the  Kronecker symbol.
\item For every positive integer $p$ let 
\begin{equation}
    \omega_p^\#:=\max_{1\leq k\leq s}\left\{\min  \left\{ | \la Q,\lambda^{(k)}\ra |;   \ Q\in {\bf N},
\, \la Q,\lambda^{(k)}\ra \ne 0,   \|Q\|< 2^p\right\}\right\}. 
\end{equation}
Then {\em Condition } $\omega^\#$ requires that 
\begin{equation}
    \sum_{p=1}^\infty \dfrac{\log({\omega^{\#}_p}^{-1})}{2^p}<\infty.
\end{equation}
\end{itemize}
\end{defin}
\begin{remark}
Introducing the set 
\begin{equation}\label{efreL}
{M}_\Lambda:=\left\{P \in {\bf N}| \  \langle P,\lambda^{(k)}\rangle=0,\quad 1\leq k\leq s
\right\},
\end{equation}
one may restate Condition AL in the form
\begin{equation}\label{eq:ALtwo}
G^{(k)}(x)=\sum_{P\in M_{\Lambda}} (x\odot G_P^{(k)}) x^P,\quad 1\leq k\leq s
\end{equation}
with each $G_P^{(k)}$ a $\mathbb C$-linear combination of the $\lambda^{(i)}$, and   $G_{\bar 0}^{(k)}=\lambda^{(k)}$.
\end{remark}
\begin{remark}\begin{itemize}
\item The series $\check v_{ik}$ are common first integrals of the $A^{(k)}$. As in the previous section one verifies that, if one normal form satisfies Condition AL, then every other normal form will.
\item Condition $\omega^\#$ is certainly satisfied whenever Condition $\omega$ holds for each individual vector field $F^{(k)}$, but is less restrictive.
\end{itemize}
\end{remark}
{The following is a version of Stolovitch's convergence theorem in our simplified scenario; compare Stolovitch, \cite{St}, Theorem 1.0.1.}
\begin{teo}\label{stolopart}
Let commuting analytic  vector fields be given as in \eqref{eq:abellie}. If some normal form satisfies Condition AL, and if Condition $\omega^\#$ holds, then there exists a convergent transformation to normal form.
\end{teo}
The proof of this theorem will mimic the proof of Theorem \ref{BSH}, with a few modifications. Our vantage point is equation \eqref{ehds2}.
We start with a nilpotency result, a variant of Lemma \ref{lem:26}.
\begin{lem}\label{lem:26L}
 If the vector fields  
 $$
 G^{(k)}(x)
 =\sum_{P\in M_\Lambda } (x \odot   G_{P}^{(k)}) x^{P}
 $$
 satisfy condition  \eqref{eq:ALone}, then $$\Delta_{G_*^{(k)}}^2=0.$$
\end{lem}
\begin{proof}
 If  condition \eqref{eq:ALone} holds, then for any term $(x \odot G_{*,P}) x^P$ of $G_*$ we see that  $G_{*,P}$ has the form \eqref{eq:ALtwo}. Observing that $P\in M_\Lambda$, Lemma \ref{lem:Deltasq} yields the assertion.
\end{proof}

\begin{lem}\label{lem_solhomolL}
Given  $\lambda^{(1)}, \dots, \lambda^{(r)}$ as before, write 
for $Q\in\bf N$, $\Vert Q\Vert\leq 2m-1$:
\[
\delta_i=\left<Q,\,\lambda^{(i)}\right>; \quad 1\leq i \leq s.
\]
If the $G^{(j)}(x)$ satisfy Condition AL, consequently
\begin{equation}\label{eq:cutgstar}
G^{(j)}_*(x)=\sum_{i=1}^s v_{ij}(x)\,\lambda^{(i)}\odot x
\end{equation}
with suitable polynomials $v_{ij}$, $v_{ij}(0)=0$, and $k$ is such that $\delta_k\not=0$, then the unique solution of equation \eqref{ehds2} is
\begin{equation}\label{eq:hdellie}
h_{\bar\delta}(x)=\dfrac1{\delta_k\left(1+\sum_i \alpha_i v_{ik}(x)\right)}\,\left({\rm Id}+  \dfrac1{\delta_k\left(1+\sum_i \alpha_i v_{ik}(x)\right)} \Delta_{G_{*}^{(k)}}\right)\,{\rm Pr}_{2m-1}(F^{(k)}_{*,\bar\delta}),
\end{equation}
where $\alpha_i:=\delta_i/\delta_k$, $1\leq i\leq s$. 
\end{lem}
\begin{proof}
By \eqref{jacmulti}, Lemma \ref{lem:hD_n} and equation \eqref{eq:cutgstar}, one has
\[
\left[G_*^{(k)},\,h_{\bar\delta}\right]=\sum_i \delta_i v_{ik}h_{\bar\delta}.
\]
Substituting this in equation \eqref{ehds2}, the assertion 
 follows with inverting the geometric series, using the nilpotency of $\Delta_{G_*^{(k]}}$.
\end{proof}
\begin{remark}\label{rem:bestk}
Choosing $\delta_k$ in the lemma of largest modulus, one has $|\alpha_i|\leq 1$ for all $i$ in \eqref{eq:hdellie}.
\end{remark}

Now we state the counterpart of Lemma \ref{lem:stepest} on estimates.

\begin{lem}\label{lem:stepestLie}
Consider the setting of Proposition \ref{pro:lien3} and Corollary \ref{cor4711}.
Assume that for some $\rho$, $\frac12<\rho\leq 1$, and $1\leq \ell\leq s$ one has 
\be \label{estlie}
 \overline { |F^{(\ell)}|}_\rho<1, \, \overline { |{G_*}^{(\ell)}|}_\rho <c_1,\,\overline {| \Delta_{ {G_*}^{(\ell)}}|}_\rho <c_1,\,
\ee
with some constant $c_1$, and moreover assume that Condition AL holds.  Furthermore let $m=2^p$. For $Q\in{\bf N}$ with $m\leq \Vert Q\Vert\leq 2m-1$ let $\delta^{(i)}=\left<Q,\,\lambda^{(i)}\right>$, $1\leq i\leq s$. 
Then there exists a constant $c_2$ such that 
$$
\overline{ |h|}_\rho <\dfrac{c_2}{{\omega^\#_{p+1}}^2}.
$$
\end{lem}
\begin{proof}
We choose $k$ as in Remark \ref{rem:bestk}. The proof runs closely parallel to the proof of Lemma \ref{lem:stepest}. 
As for part (i) of that proof, we have the estimate 
\[
| \delta^{(i)}/\delta^{(k)}|\leq 1,
\]
an analog to \eqref{eq:c_est}. Then one imitates the arguments in the proof of Lemma \ref{lem:stepest} to 
show that there is a $c_1$ such 
\be \label{eq:asrLie}
 \overline{|v_{ik}|}_\rho < \frac{1}{2r }.
\ee
The proof of part (ii) carries over with obvious modifications in the notation, except for the crucial involvement of $\omega_{p+1}^\#$, due to our choice of $k$.
\end{proof}
From here on, one may (as in the proof of Theorem \ref{BSH}) appeal to Bruno \cite{BruB}: We have shown (modulo notation) that the statement of the Corollary on p.\ 204 of \cite{BruB} holds. Again, the remainder of the proof is identical (modulo notation) to the proof of \cite{BruB}, Theorem 4 on p.\ 204 ff.

{
\begin{remark}
    Our proof of convergence in case of Bruno's $A_2$ is different from the one given in Stolovitch \cite{St}, who invokes the existence of a nontrivial analytic commuting vector field.
\end{remark}}

\section{About Bruno's Case 1 }\label{sec:Bruno1}

Bruno's Case 2 (Condition $A_2$ as well as the Simplified Condition A) does not cover the classical convergence result due to Poincar\'e, among others. We will deal with some of these now, mostly following the pattern in Bruno \cite{BruB}.

For Bruno's Case 1, which is complementary to Case 2, we will assume that there is an $l\geq 0$ such that $\lambda_1, \dots , \lambda_l$  belong to a
real line $L=\mathbb R\omega$ passing through the origin of $\C$  whereas $\lambda_{l+1}, \dots, \lambda_n$ 
 lie on the same side of $L$. (For $l=0$ all eigenvalues lie to one side of $L$; this is the case considered by Poincar\'e.)

 Because the notions and computations regarding normal forms, as well as convergence results, remain unaffected upon multiplication of all eigenvalues by a nonzero complex constant, we may -- and from here on will -- assume that $L=i\,\mathbb R$ is the imaginary line, and moreover that $\lambda_{l+1}, \dots, \lambda_n$  have positive real parts. (In dynamical systems language, the eigenspaces for the latter eigenvalues span the unstable manifold of the system.)

Then we can write 
$$
 \lambda=\mu+i \nu,
$$
where $\mu$ and $\nu$ are real vectors, and 
\be \label{eq:6}
0=\mu_1=\dots =\mu_l< \mu_{l+1}\le \dots \le \mu_n
\ee 
with appropriate numbering. We furthermore abbreviate $\mu_*:=\mu_{l+1}$, $\mu^*:=\mu_n$ for the minimum and maximum of the nonzero $\mu_j$.

From here on, for any $n$-tuple $v=(v_1, \dots, v_n ) $ we denote by 
$v'$ the $l$-tuple $(v_1, \dots, v_l ) $ and by 
$v''$ the $n-l$-tuple $(v_{l+1}, \dots, v_n )$. We will also write $v=\begin{pmatrix}
    v'\\ v''
\end{pmatrix}.$ 
The notation $v\ge 0$ means that $v_i \ge 0 $
for all $i=1,\dots, n$.

For real $Q$ the equation 
$$
\la Q, \lambda \ra =0
$$
is equivalent to the system 
\be \label{eq:8}
\la Q, \mu \ra =0,
\ee
$$
\la Q, \nu \ra =0.
$$

By \eqref{eq:6} we have $\mu'=0$ and $\mu''>0$
and  \eqref{eq:8} is equivalent to the equation 
\be \label{eq:11}
\la Q'', \mu''\ra =0.
\ee
We state and prove some facts regarding condition \eqref{eq:11}. Denote by ${\bf N}''$ the restriction of $\bf N$ to $\mathbb R^{n-l}$.
\begin{lem}\label{lem:poinc}
For all $Q\in {\bf N}$: 
    \begin{enumerate}[(a)]
        \item There are positive constants $c_1,\,c_2$ such that
        \[
        |\la Q'',\,\mu''\ra| \geq c_1|Q''|\quad\text{whenever}\quad |Q''|\geq c_2.
        \]
        \item There are only finitely many $Q''\in {\bf N}''$ such that $\la Q'',\,\mu''\ra=0$.
        \item There exists a $c_3>0$ such that 
        \[
        |\la Q'',\,\mu''\ra|\geq c_3|Q''|\quad\text{whenever}\quad \la Q'',\,\mu''\ra\not=0.
        \]
    \end{enumerate}
\end{lem}
\begin{proof}
    For part (a) we consider nonnegative integers $p_{l+1},\ldots, p_n$ and set $Q''=(p_{l+1},\ldots,p_n) -e_k$ with some $k>l$. Comparing $\sum p_j\mu_j-\mu_k$ to $\sum p_j+1$, we have 
    \[
    \dfrac{1}{\mu_*}(\sum p_j\mu_j-\mu_k)\geq \sum p_j-\dfrac{\mu^*}{\mu_*}.
    \]
    Moreover 
    \[
    \sum p_j-\dfrac{\mu^*}{\mu_*}\geq \frac12\,(\sum p_j+1)\quad{if}\quad \sum p_j+1\geq 3+2\dfrac{\mu^*}{\mu_*}=:c_2.
    \]
    The assertion follows with $c_1=\mu_*/2$.\\
    Parts (b) and (c) follow directly because the set of all $Q''$ with $|Q''|< c_2$ is finite.
\end{proof}


We take a closer look at the condition $\la Q,\,\mu\ra=0$ for $\|Q\|=0$. Setting $Q=(p_1,\ldots,p_n)-e_k$, the condition implies:
\begin{enumerate}[(i)]
    \item When $k\leq l$ then $p_{l+1}=\cdots=p_n=0$;
    \item When $k>l$ then $p_j>0$ for some $j>l$.
\end{enumerate}
   
Lemma \ref{lem:poinc} implies that any solution  $Q=(Q',Q'')\in {\bf N}$ to \eqref{eq:8} belongs 
to one of the following sets:
\be \label{M123}
\begin{aligned}
  {M}_1^*= &\{Q: Q''=0  \},\\
{M}_2^*= &\{Q: Q''=e_j-e_k,  l< k,j\le n  \},\\ 
{M}_3^*= &\{Q: Q''=\sum_{j=l+1}^m p_j e_j -e_k, \  {\rm where \ } p_j\ge 0, m<k, \mu_m <\mu_i   \},  
\end{aligned}
\ee 
that is, 
$$
{M}_\lambda={M}_1^*\cup {M}_2^*\cup {M}_3^*.
$$  
 Note that 
  \be \label{QM23}
 \| Q''\|=0  \quad\text{for} \quad Q\in {M}_2^*; \quad
 \|Q''\|>0\quad \text{for} \quad Q\in {M}_3^*.
 \ee
 Moreover,  a term $e_j-e_k$ appears in $M_2^*$ only if $\mu_j=\mu_k$. 
Then, in view of (i) above, the normal form of system \eqref{sys} can be written as 
\be \label{nf_c1simple}
\dot y=
 \sum_{Q\in {M}_1^*} {y' \ch y''} \odot {G'_Q\ch G''_Q} y^Q +\sum_{Q\in {M}_2^*\cup M_2^*} {y'\ch y''}\odot {0 \ch G''_Q} y^Q. 
\ee
{
A refined version of \eqref{nf_c1simple} can be written as 
\be \label{nf_c1}
\dot y=
 \sum_{Q\in {M}_1^*} {y' \ch y''} \odot {G'_Q\ch G''_Q} y^Q +\sum_{Q\in {M}_2^*} y\odot {E}_{Q} y^{Q}+   \sum_{Q\in {M}_3^*} {y'\ch y''}\odot {0 \ch G''_Q} y^Q. 
\ee
Note that for $Q\in {M}_2^*$  one has ${E}_{Q} y^{Q}= e_i y^{Q'} y_j y_i^{-1}$ ($l\le i,j$) and $E_Q=0$ if $\mu_i\ne \mu_j$. }
{
\begin{remark}\label{rem:echelon}
  From (i) and (ii) above and \eqref{nf_c1} one sees that the normal form system is triangular; i.e.,
\be\label{nfBap}
\begin{array}{rcl}
    \frac{d}{dt}y'&=& G'(y')   \\
    \frac{d}{dt} y''& =& G''(y',y'')=B(y')y''+\widehat{B}(y',y'')
\end{array}
\ee
with formal power series $G',\, G''$ satisfying $G'(0)=0$ and $G''(y',0)=0$, $ \widehat B(y) $  contains neither constant nor linear  terms in  $y''$, and 
\begin{equation}\label{eq:Bij}
    B(y')=\{ b_{ij}(y') \}_{l+1 \le i,j\le n }
\end{equation}
corresponding to the terms from $M_1^*$ and $M_2^*$. Moreover, $B$ is in ``block diagonal'' form, with  $b_{ij}=0$ whenever $\mu_i\not=\mu_j$. This implies for the analytic case that $y'=0$ defines an invariant manifold (the unstable manifold), and also that $y''=0$ defines an invariant manifold (center manifold) of the system. 
\end{remark}
}


Following Bruno \cite{BruB}, there are two subcases to consider, denoted as 1* and 1**:

1*: $\lambda_1,\dots, \lambda_l$   are pairwise commensurable;

1**: there is an incommensurable pair  among $\lambda_1,\dots, \lambda_l$. \\

{We will discuss Case 1* (which includes the Poincar\'e case, and gives an impression of the general proof strategy) in detail, but we will not delve into Case 1**, since our approach does not lead to substantial simplification here. (Bruno already used the Lie algebra structure in treating this case.) Also, we will briefly make notice of the natural invariant manifolds appearing in Case 1.}


We write \eqref{nf_c1simple} in the form 
\begin{multline} \label{nf_c1*p}
\dot x= 
G(x)={G'(x) \ch G''(x)}= x \odot  \lambda+ 
\sum_{P\in {M}_1^*, P\ne 0} {x' \ch x''} \odot {G'_P\ch 0} x^P +\\  \sum_{P\in   {M}_2^* \cup {M}_3^*} {x'\ch x''}\odot {0 \ch G''_P} x^P = x\odot \lambda +G_\#(x).
\end{multline}
\begin{teo}\label{Thm:Bruno1*}
If $\lambda$ satisfies condition  1* and  
a normal form of the analytic system \eqref{sys} satisfies 

\be\label{eq:cA}
 G'(x') = (\lambda' \odot x') a(x'),
\ee
where $a(x')=1+\dots $ is a power series,
then it admits a convergent 
transformation into normal form.
\end{teo}

To start the proof, we turn to the homological equation. Let  ${F_*}$ be as in Proposition \ref{pro:1n},
$h$ be as in Corollary \ref{cor47}, and $\delta$ a nonzero eigenvalue of ${\rm ad}_A$.
By \eqref{ehds} we have
\be \label{eq92}
\delta  h_\delta - \prtm ( [h_\delta, G_*])+\prtm ( F_{*,\delta})=0\quad \text {in}\,\,\mathcal V^{2m-1},
\ee
where $G_*$ consists of terms of $G_\#$ defined by \eqref{nf_c1*p} of order at most $m-1$. 
{As for the plan of the proof, we note that Remark \ref{hdhat} applies here as well.}

Using 1) of Lemma \ref{lem:hD_n} and \eqref{DelU} we write the latter equality as 
\be \label{eq:hd}
-\delta h_\delta -\Delta_{ h_\delta} G_* +\Delta_{G_*} h_\delta - F_\delta =0,
\ee
where 
$
 F_\delta=\prtm (F_{*,\delta}).
$

\begin{lem}\label{lem:onestar}
\begin{enumerate}[(a)]
    \item There exists an $\epsilon>0$ such that $|\delta|>\epsilon$ for every nonzero eigenvalue of ${\rm ad}_A$. In particular, Condition $\omega$ is satisfied.
    \item Moreover, there exists $c>0$ such that
$|Q''|<c |\delta|$ for all $Q''\in {\bf N}''$.
\end{enumerate}
\end{lem}
\begin{proof}
We have 
\[
\delta =\la Q,\lambda\ra=\la Q'',\mu''\ra + i\la Q',\nu'\ra +i\la Q'',\nu''\ra.
\]
Therefore, with Lemma \ref{lem:poinc} there exists a $\epsilon_1 >0$ such that
\[
|\delta|=|\la Q,\lambda\ra|\geq |\la Q'', \mu''\ra|\geq \epsilon_1\,|Q''|\geq \epsilon_1\quad\text{if  }\la Q'',\mu''\ra\not=0.
\]
By the same lemma, there exist only finitely many $Q''$ such that $\la Q'',\mu''\ra=0$, for which one gets
\[
\la Q,\lambda\ra =i(\la Q',\nu'\ra+\la Q'',\nu''\ra).
\]
We show that for these $Q=(Q',Q'')$ there exists some $\epsilon_2>0$ such that $|\la Q,\lambda\ra|\geq \epsilon_2$ unless $\la Q,\lambda\ra=0$. By the defining condition for case 1*, there exists a real $\beta >0$ such that $\widetilde \nu:=\beta \nu'$ has integer components. Now fix some $Q''$ and set $\kappa:=\la Q'',\mu''\ra$. One sees that
\[
\la Q',\nu'\ra+\kappa \in \beta\left(\mathbb Z+ \kappa/\beta\right),
\]
and this set is discrete. Since there are only finitely many $Q''$, the set of all $\la Q,\lambda\ra =i(\la Q',\nu'\ra+\la Q'',\nu''\ra) $, with $\la Q'',\mu''\ra =0$ is discrete, hence there exists a lower bound $\epsilon_2>0$ for the absolute values of its nonzero elements.\\
Altogether one has $|\la Q,\lambda\ra|\geq {\rm min}\,\{\epsilon_1,\,\epsilon_2\}$ unless  $\la Q,\lambda\ra=0$, and part (a) is proven.\\
For part (b), the case $\la Q'',\mu''\ra\not=0$ is taken care of. When $\la Q'',\mu''\ra=0$ but $\la Q,\lambda\ra \not=0$ then $\epsilon_2/({\rm max}\,\{|Q''|\}$ satisfies the asserted inequality.
\end{proof}

\begin{lem}\label{lem:cA}
Assume that for system \eqref{nfm} 
\be \label{est1}
\frac 12 <\rho\le 1,  \ \overline{|F|}_\rho<1, \ \overline{| {G_*}|}_\rho <c_1,   \ 
\overline{ |\Delta_{G_*}|}_\rho <c_1
\ee
and its normal form satisfies \eqref{eq:cA},
then there is a constant $c_2$, independent of $m$, such that
\be \label{est_h}
\overline{| h|}_\rho <c_2.
\ee
\end{lem}
\begin{proof}
In view of \eqref{eq:cA} we have 
$$
G'_P = \lambda' a_P \qquad (P\in M_1^*) .  
$$
Therefore, we can write the projection of $G(x)$ from  \eqref{nf_c1*p} to ${\mathcal V}^{m}$ as
\begin{multline*}
\prm\left(  {0\ch \lambda''}\odot {0\ch x''} + \sum_{P\in {M}_1^*} (x \odot \lambda) a_P x^P  -\right. \\ \left. \sum_{P\in {M}_1^*} {x' \ch x''} \odot {0 \ch  \lambda''} a_P x^P + \sum_{P\in 
{M}_2^*\cup {M}_3^*} {x'\ch x''}\odot {0 \ch G''_P} x^P\right) =\\
\prm \left(x\odot \lambda+ \sum_{P\in {M}_1^*, P\ne \bar 0} (x \odot \lambda) a_P x^P + \sum_{P\in  
{M}_\lambda} {x'\ch x''}\odot {0 \ch \breve {G}''_P} x^P\right) =\\ x\odot \lambda +G_1(x) +G_2(x),
\end{multline*}
where 
$$
 \breve {G}''_P= G''_P -\lambda'' a_P
$$
and 
\be \label{G1G2*}
G_1(x)=\prm\left(\sum_{P\in {M}_1^*, P\ne \bar 0} (x \odot \lambda) a_P x^P\right), \  G_2(x)=\prm\left( \sum_{P\in  
{M}_\lambda, P\ne \bar 0} {x'\ch x''}\odot {0 \ch \breve {G}''_P} x^P\right).
\ee
Now let $\delta$ be a nonzero eigenvalue of ${\rm ad}_A$, and consider all $Q$ of order greater or equal to  $m$ and  smaller or equal to  $ 2m-1$  with $\la Q,\lambda\ra=\delta$; thus we write $h_\delta$  appearing in \eqref{eq:hd} as 
\[
h_\delta=\sum_Q(x\odot H_Q)x^Q.
\]
Using 2) of  Lemma \ref{lem:hD_n}
we have 
\begin{multline} \label{eq:c1*}
[h_\delta, G_*]= [h_\delta, G_1+G_2]= \\
  \Delta_{G_1}\hd-
\Delta_{\hd} G_1+\Delta_{G_2} \hd -\Delta_{\hd} G_2 =\\ 
\Delta_{G_*}  h_\delta- \delta\hd  \tilde a(x) - \Delta_{\hd} G_2,
\end{multline}
with $\tilde a(x)=a(x)-1$. 
Since 
\be\label{hdG2}
\Delta_{\hd} G_2 = \sum_{P,Q} \la Q'',\breve{G_P}''\ra (x\odot H_Q) x^{P+Q}, 
\ee
 from  \eqref{eq:hd} and \eqref{eq:c1*} we obtain 
\be \label{eq:hd1}
 h_\delta=  - h_\delta \tilde a - \frac{1}{\delta} \sum_{P,Q} \la Q'',\breve{G_P}''\ra (x\odot H_Q) x^{P+Q}  +\frac{1}{\delta} \Delta_{G_*} h_\delta + \frac{1}{\delta} F_\delta
\ee
where {\it  in \eqref{hdG2}, \eqref{eq:hd1} and below the summation is over $P$ and $Q$ such that $\|P+Q\|\le 2m-1.$} 
By Lemma \ref{lem:onestar} there exist $\epsilon>0$ such that $|\delta |>\epsilon$, and $c>0$ such that $|Q''|<c|\delta|$ (both constants independent of $m$).

Therefore, 
\begin{multline}\label{hd1*}
\sum_{P,Q} |\la Q'',\breve{G_P}''\ra| (x\odot |H|_Q) x^{P+Q} \prec \sum_{P,Q} | Q''| |\breve{G_P}''| (x\odot |H|_Q) x^{P+Q}\prec \\
 \sum_{P,Q} c | \delta | |\breve{G_P}''| (x\odot |H|_Q) x^{P+Q}\prec \sum_{P,Q} c | \delta | 
 ( (x\odot |H|_Q) x^{Q} ) (|  \lambda'' a_P-G''_P  | x^P).
\end{multline}
where $ |H|_Q=(|H_Q^{1}|, \dots |H_Q^{n}| )^\top $.

Summing up over $\delta$ and taking into account \eqref{hd1*} and the observation  that $|\delta|^{-1} < \epsilon^{-1}$   we have 
\be \label{suh}
  |\bar h| \prec  |\bar h| |\bar a|+ c |\bar h | (|\lambda| |\bar a|+ |G''| )  + \frac{1}{\epsilon}|\overline{\Delta_{G_*}}| |\bar h| +\frac{1}{\epsilon} |\bar F|.
\ee

Since $\rho >1/2$ from \eqref{est1} we have 
$$
|\bar G|_\rho < 2 c_1, \quad |\lambda'| |\bar a|_\rho <2 c_1.
$$
If $\lambda'=0$ we set $\hat \lambda =1, \ a=0$;
otherwise we set
$\hat{\lambda}= |\lambda|.
$
Then, from \eqref{suh} we obtain 
$$
  |\bar h|_\rho \prec  |\bar h|_\rho |\bar a|+ c |\bar h |_\rho ( 2 c_1 +|\lambda| c_1 \hat \lambda^{-1} )  + \frac{1}{\epsilon}|\overline{\Delta_{G_*}}|_\rho  |\bar h|_\rho +\frac{1}{\epsilon} |\bar F|_\rho.
$$

Taking into account \eqref{est1} we have 
$$
  |\bar h|_\rho \prec  |\bar h|_\rho   \left(2 c_1 \hat\lambda^{-1} + c( 2 c_1 |\lambda| c_1 \hat\lambda^{-1}+ 2 \frac{c_1}{\epsilon} \right) + 
\frac 1\epsilon |\bar F|_\rho.
$$
Taking $c_1$ such that 
$$ c_1 \left(2 c \hat\lambda^{-1} +  2 c_1 |\lambda| c_1 \hat\lambda^{-1}+ 2 \frac{1}{\epsilon} \right)  = \frac{1}{2} $$
we see that \eqref{est_h} holds with $c_2=2/\epsilon$.
\end{proof}

To finish the proof of the theorem, repeat the arguments at the end of subsection \ref{A2est}.


\begin{remark}
     Bruno's Case 1* also includes the borderline case with $l=0$; thus $x=x''$ and all eigenvalues lie in the open positive half plane, which is the case first considered by Poincar\'e. Here Theorem \ref{Thm:Bruno1*} holds without any conditions on the normal form. The above proof only requires obvious modifications (for instance, $M_1^*$ is the empty set, and the summations attain a simpler form).
\end{remark}
Bruno's Case 1** deals with the scenario when $\lambda_1,\ldots,\lambda_l$ are not pairwise commensurable. To ensure convergence, in addition to condition $\omega$ for these eigenvalues, Bruno imposes a very restrictive algebraic condition, called Condition $A_1''$ , on  the matrix $B$ in the representation \eqref{nfBap}. (In our framework, Condition
\(A_1''\) can be seen as a condition closely related to Stolovitch's notion of good deformation; see Definition 6.3.10 of
\cite{St}.)

The fact that convergence to full normal form is often not guaranteed in Case 1** is mitigated by the fact that convergence to a partial normal form (which is often sufficient to understand the dynamics of the system) can often be ensured.

Given Bruno's Case 1 (Cases 1* as well as 1**), any analytic system \eqref{Asn} admits a convergent transformation to normal form on the (invariant) unstable manifold. In detail: 
\begin{itemize} \item There exists a convergent transformation to a form
\be
\begin{array}{rcl}\label{nfim}
    \frac{d}{dt}z'&=& \widehat F'(z',z'')   \\
    \frac{d}{dt} z''& =& \widehat F''(z',z'')
\end{array}
\ee
such that $F'(0,z'')=0$ (thus the subset given by $z'=0$ is invariant), and the differential equation 
\[
 \frac{d}{dt} z''=\widehat F''(0,z'')
\]
is in normal form.
\item Moreover, if the eigenvalues $\lambda_1,\ldots,\lambda_l$ satisfy the Siegel-Pliss condition, and for some formal normal form \eqref{nfBap} one has $G'(y')=a'(y')\lambda'\odot y'$ with a power series $a'$ (i.e., simplified condition A holds for $G'$), then there is a convergent transformation to a form \eqref{nfim} with $F''(z',0)=0$, so $z''=0$ defines an invariant set.
\end{itemize}
Both these statements can be proven as a consequence of Bibikov's theorem about normal forms on invariant manifolds; see \cite{Bibikov}, Theorems 3.1 and 3.2.

Finally, we remark that in Chapter III of his monograph \cite{Bru}, Bruno introduces ``sets of analyticity'' in subsection 3.3 ff.; this also leads to certain invariant sets.

\section{Local and formal integrability}\label{sec:integrable}

In this section we consider a single vector field \eqref{sys}, with the linear part $A$ in Jordan form.  Let $A_s$ be the semisimple part of $A$, so $A_s x= \lambda \odot x. $ 
{Throughout this section we will assume that $\lambda\not=0$.}
 We recall
\[
M_\lambda^+= \left\{ P\in \mathbb Z^n_{\ge 0};\,\la P,\lambda\ra=0\right\}, 
\]
and let $R_\lambda$ be the submodule spanned by  $M_\lambda^+$. Moreover we introduce
\begin{equation}\label{barrlambda}
\bar R_\lambda= \left\{ P\in \mathbb Z^n;\,\la P,\lambda\ra=0\right\}.
\end{equation}
Both $R_\lambda$ and $\bar R_\lambda$ are free $\mathbb Z$-modules, being torsion-free submodules of $\mathbb Z^n$, and $M_\lambda^+$ is an abelian monoid whose  elements correspond to monomial first integrals of $\dot x =A_s x.$ 

We will characterize formal, and formally meromorphic, first integrals of normal forms in this section, and also discuss their relevance for convergence.
Thus let 
\be \label{nfB1}
\dot x=  G(x),\quad \text{where  }G(x)=\sum_{P\in {M}_\lambda}  (x\odot G_{P}) x^{P}
\ee
be in normal form. We recall:
\begin{lem}\label{lem:it}
For $Q\in \Z^n$, $x^Q$ is a first integral of \eqref{nfB1} if and only if 
\be \label{prPG}
\la Q, G_P \ra =0 \quad \text{for all  } P\in M_\lambda
\ee
\end{lem}
\begin{proof}
Taking the Lie derivative of $x^Q$ with respect to the vector field \eqref{nfB1}
we have 
$$
L_G(x^Q)=   x^Q \left(  \sum_{P\in {M}_\lambda}  \la Q, G_P \ra  x^P \right).  
$$
\end{proof}

\subsection{Generalities}

We denote by $\mathbb C[[x_1,\ldots,x_n]]$ the ring of formal power series. See Lang \cite{Lang} Ch.\ IV \S9 or Zariski and Samuel  \cite{ZarSam}, Ch.~VII for its properties; notably this is a unique factorization domain. By $\mathbb C((x_1,\ldots,x_n))$ we denote its field of fractions; we will call its elements {\em formally meromorphic functions}, as in Llibre et al.\ \cite{LWZ}. (Other names appear in the literature; for instance Cong et al. \cite{CLZ} speak of generalized rational functions.) 

We first recall and slightly extend some results about formal first integrals of normal forms; see Walcher \cite{WalcherNF}, Zhang \cite{XZ}, Llibre et al.\ \cite{LPW,LWZ}; most of these are essentially known. As a preliminary step we characterize first integrals of the semisimple linear part.
\begin{pro}\label{formfi}
\begin{enumerate}[(a)]
    \item  Every formal first integral of $\dot x=A_sx$ can be represented as a series in monomials $x^P$ with $P\in M_\lambda^+$. The number of independent formal first integrals is equal to the rank of $R_\lambda$.
    \item Every formally meromorphic first integral of $\dot x=A_sx$ can be represented in the form $\phi(S(x))$, with $\phi$ a formal power series {in $d$ variables} and
    \[
    S(x)=\begin{pmatrix}
        x^{S_1}\\ \vdots\\ x^{S_d}
    \end{pmatrix};
    \]
    with suitable $S_1,\ldots,S_d\in \bar R_\lambda$. Here $d$ is the rank of $\bar R_\lambda$.
\end{enumerate}
\end{pro}
\begin{proof}
For part (a) see Llibre et al \cite{LPW}. For part (b) see Llibre et al. \cite{LWZ}, part (iii) in the proof of Theorem 5, and also Cong et al. \cite{CLZ}.
\end{proof}
First we consider formal power series.
\begin{pro}\label{pro:LPW}
 Let the formal vector field  \eqref{nfB1} be in Poincar\'e-Dulac
  normal form. Then the following hold. 
  \begin{enumerate}[(a)]
 \item Every formal first integral of \eqref{nfB1} is also a first integral of the linear vector field
   \be\label{Als} \dot  x= A_sx= x  \odot \lambda. \ee 
 \item Assume that \eqref{nfB1} admits $d$ independent formal first integrals. Then \eqref{Als} admits $d$ independent polynomial first integrals.
 \item When ${\rm rank\,}R_\lambda=d$, then \eqref{nfB1} admits $d$ independent formal first
integrals if and only if either of the following holds: 
\begin{itemize}
    \item The differential  equation \eqref{nfB1} in normal form admits $d$ independent polynomial first integrals of system \eqref{Als};
    \item the differential equation \eqref{nfB1} in normal form admits every polynomial first integral of system \eqref{Als}.
\end{itemize} 
 \end{enumerate}
 \end{pro}
 \begin{proof}
     For part (a) see \cite{WalcherNF}, Prop.~1.8. 
     {For the sake of completeness we include a detailed proof of part (b), which seems not readily available in the literature.}  Thus assume that $\dot x=G(x)$ admits $d$ independent formal power series first integrals
     \[
     s_i=\sum_{k=0}^\infty s_{ik},\quad 1\leq i\leq d.
     \]
     By definition of (functional) independence the Jacobian
     \[
     \left(\dfrac{\partial{s_i}}{\partial x_j}\right)_{1\leq i\leq d, \,1\leq j\leq n}
     \]
     has full rank $d$ as a matrix over $\mathbb C((x_1,\ldots,x_n))$. Equivalently (with no loss of generality up to a permutation of indices) 
      we have that
     \[
     \Theta:=\det\left(\dfrac{\partial{s_i}}{\partial x_j}\right)_{1\leq i,j\leq d}\in \mathbb C((x_1,\ldots,x_n))\setminus \{0\}.
     \]
     With 
     \[
        s_i^{(m)}=\sum_{k=0}^m s_{ik},\quad \Theta_m = \det\left(\dfrac{\partial{s_i^{(m)}}}{\partial x_j}\right)_{1\leq i,j\leq d}
     \]
     we have $s_i=s_i^{(m)}+\widetilde s_i$, with all terms in $\widetilde s_i$ having order $>m$. Moreover 
     \[
     \Theta=\det\begin{pmatrix} \dfrac{\partial{s_1^{(m)}}}{\partial x_1}&\cdots& \dfrac{\partial{s_1^{(m)}}}{\partial x_d}\\
         \dfrac{\partial{s_2}}{\partial x_1}&\ldots &\dfrac{\partial{s_2}}{\partial x_d}\\
         \vdots &  & \vdots \\
         \dfrac{\partial{s_d}}{\partial x_1}& \ldots & \dfrac{\partial{s_d}}{\partial x_d}
     \end{pmatrix}
     + \det\begin{pmatrix} \dfrac{\partial \widetilde s_1}{\partial x_1}&\cdots& \dfrac{\partial \widetilde s_1}{\partial x_d}\\
         \dfrac{\partial{s_2}}{\partial x_1}&\ldots &\dfrac{\partial{s_2}}{\partial x_d}\\
         \vdots &  & \vdots \\
         \dfrac{\partial{s_d}}{\partial x_1}& \ldots & \dfrac{\partial{s_d}}{\partial x_d}
     \end{pmatrix}
     \]
     by linearity of the determinant in the first row. Now all terms appearing in the first row of the second matrix are of order $>m-1$; expansion along the first row shows that all terms appearing in the second determinant have order $>m-1$. Briefly,
       \[
     \Theta=\det\begin{pmatrix} \dfrac{\partial{s_1^{(m)}}}{\partial x_1}&\cdots& \dfrac{\partial{s_1^{(m)}}}{\partial x_d}\\
         \dfrac{\partial{s_2}}{\partial x_1}&\ldots &\dfrac{\partial{s_2}}{\partial x_d}\\
         \vdots &  & \vdots \\
         \dfrac{\partial{s_d}}{\partial x_1}& \ldots & \dfrac{\partial{s_d}}{\partial x_d}
     \end{pmatrix}
     + o(m-1).
     \]
     Repeating this argument for rows $2,\ldots,d$ yields
     \[
     \Theta=\Theta_m+ o(m-1).
     \]
     Now let $m$ be such that the lowest term in $\Theta$ has order $\leq m-1$. Then the lowest order term of $\Theta_m$ must cancel the lowest order term of $\Theta$, and therefore $\Theta_m\not=0$. This means that the $s_i^{(m)}$ are independent polynomial first integrals of $\dot x=Ax$.\\
    For part (c) see Llibre et al.\ \cite{LPW}, Proposition 5.
 \end{proof}
 
In view of Propositions \ref{formfi} and \ref{pro:LPW} we say that system \eqref{sys}  is {\it  maximal formally integrable} if ${\rm rank\,}R_\lambda=d$ 
and \eqref{sys} admits  $d$ independent formal 
first integrals. (In Llibre et al.\ \cite{LPW} this scenario was referred to as ``conservation of all first integrals''.)
\begin{pro} \label{pro:69}
Let  system \eqref{sys} be given, with  \eqref{nfB1} a corresponding normal form. 
\begin{enumerate}[(a)]
    \item  Let   ${\rm rank\,}R_\lambda=d$ . Then 
    the  system is maximal formally integrable if and only if
\be \label{eq:qgp}
\la Q, G_P \ra =0 \qquad \text{for all  } P\in M_\lambda, \,\text{all  }  Q\in M_\lambda^+.
\ee
\item Let the system be maximal formally integrable, and assume for the matrix $A$ that $R_\lambda $ has rank $d$, with basis $Q_1,\ldots,Q_d$. Then there exist linearly independent
\[
\lambda^{(1)},\,\lambda^{(2)},\ldots,\lambda^{(n-d)}\in \mathbb C^n
\]
such that all $\la Q_i,\,\lambda^{(j)}\ra=0$ and the normal form satisfies 
  \be \label{Gdec}
    G(x)=\sum_{j=1}^{n-d} (\gamma_j+\check s_j(x))\cdot \lambda^{(j)}\odot x
    \ee   
    with constants $\gamma_1,\ldots,\gamma_{n-d}$ 
    and 
    $$
    \check s_j(x)=\sum_{P\in M_\lambda }  \sigma_{j,P} x^P, \quad   s_j(0)=0, \quad j=1,\dots, n-d.  
    $$
\end{enumerate}
\end{pro}
\begin{proof} 
We first show the nontrivial direction of part (a). By Proposition \ref{pro:LPW}, since system \eqref{nfB1} admits $d$ formal 
first integrals, it also admits every monomial first integral of \eqref{Als}.
Therefore, 
$x^Q$ is a first integral of \eqref{nfB1} for any $Q\in M_\lambda^+$. Using Lemma \ref{lem:it} we see that \eqref{eq:qgp} holds. 

Turning to part (b), consider the $\mathbb C$-vector space of all $\mu\in \mathbb C^n$ such that 
\be \label{eq:Ud}
\la\mu,Q_1\ra=\cdots=\la\mu,Q_d\ra=0.
\ee
This space is the kernel of the linear map from $\mathbb C^n$ to $\mathbb C^d$ whose matrix has rows $Q_1,\ldots,Q_d$ (thus is
of rank $d$), and therefore has dimension $n-d$. Moreover, it
contains $\lambda$ by construction.

Let 
\[
\lambda=\lambda^{(1)},\,\lambda^{(2)}\ldots,\lambda^{(n-d)}
\]
be a basis. By part (a) every term $(x\odot G_P) \,x^P$ in the normal form satisfies $\la Q_j,G_P\ra=0$ for $1\leq j\leq d$, whence $G_P$ is a linear combination of $\lambda_1,\ldots,\lambda_{n-d}$, and this term can be written as $\sum (\sigma_{i,P}\lambda^{(i)}\odot x) x^P)   $, so the corresponding term in $\check s_j $ is $\sigma_{j,P} x^P $.
\end{proof} 
 {Now we turn to formally meromorphic first integrals. First we extend the notion of isoresonance (compare Definition \ref{def:lambdadecomp}).
\begin{defin}
\begin{enumerate}[(i)]
    \item We say that the {set $\{\lambda^{(1)},\ldots,\lambda^{(r)}\}$ } is {\em rational isoresonant} with respect to $\lambda$ if every rational first integral of $Ax$
        is also a first integral of all $A_jx$, $1\leq j\leq r$. 
    \item {We say that a vector field \eqref{nfB1} in normal form, with an additive decomposition \eqref{ladec}, satisfies the {\em rational isoresonance condition} if $\{\lambda^{(1)},\ldots,\lambda^{(r)}\}$ is rational isoresonant with respect to $\lambda$.} 
\end{enumerate}
        \end{defin}
}
{
\begin{remark}
    By a special instance of Proposition \ref{formfi} (b), one obtains an equivalent characterization: $\{\lambda^{(1)},\ldots,\lambda^{(r)}\}$  is {\em rational isoresonant} with respect to $\lambda$ if and only if 
    \be\label{condFMM} \left<P,\lambda\right>=0\ \Longrightarrow\   \left<P,\lambda^{(j)}\right>=0 \quad \forall  \ P\in \mathbb Z^n.
    \ee
    There is yet another equivalent characteriztion: The rational isoresonance condition holds if and only if every Laurent monomial first integral of $Ax$
        is also a first integral of all $A_jx$, $1\leq j\leq r$. 
\end{remark}
}

\begin{pro}\label{pro:LWZ}
 Let the formal vector field  \eqref{nfB1} be in Poincar\'e-Dulac
  normal form. Then the following hold. 
  \begin{enumerate}[(a)]
 \item Every formally meromorphic first integral of \eqref{nfB1} is also a first integral of the linear vector field \eqref{Als}.
 \item Assume that \eqref{nfB1} admits $d$ independent formally meromorphic first integrals. Then \eqref{Als} admits $d$ independent Laurent monomial first integrals.
 \item When ${\rm rank\,}\bar R_\lambda=d$, then \eqref{nfB1} admits $d$ independent formally meromorphic first
integrals if and only if either of the following holds:
\begin{itemize}
    \item The differential  equation \eqref{nfB1} in normal form admits $d$ independent Laurent monomial first integrals of system \eqref{Als};
    \item the differential equation \eqref{nfB1} in normal form admits every Laurent monomial first integral of system \eqref{Als}.
    \end{itemize} 
    {In this scenario, we call the system {\em maximal formally meromorphic integrable}.}
    \item  
    {{Assume that for the given matrix $A$, $\bar R_\lambda $ has rank $d$, with basis $Q_1,\ldots,Q_{d}$, and that the system admits $d$ independent formally meromorphic first integrals. Then there exist linearly independent
\[
\lambda^{(1)},\ldots,\lambda^{(n-d)}\in \mathbb C^n
\]
such that all $\la Q_i,\,\lambda^{(j)}\ra=0$ and the normal form satisfies 
  \be \label{Gdecrat}
    G(x)=\sum_{j=1}^{n-d} (\gamma_j+\check s_j(x))\cdot \lambda^{(j)}\odot x.
    \ee   
Thus $\lambda=\gamma_1\lambda^{(1)}+\cdots+\gamma_{n-d}\lambda^{(n-d)}$ is an additive decomposition that satisfies the rational isoresonance condition.
    Moreover}
    $$
    \check s_j(x)=\sum_{P\in M_\lambda }  \sigma_{j,P} x^P, \quad   \check s_j(0)=0, \quad j=1,\dots, n-d.  
    $$}
 \end{enumerate}
 \end{pro}
\begin{proof}
    For part (a), see Llibre et al. \cite{LWZ}, Theorem 1. For part (b), we use and expand the arguments in \cite{LWZ}, part (iii) in the proof of Theorem 5:  The vector $R$ of independent meromorphic first integrals of the normal form, due to part (a), consists of formally meromorphic first integrals of the semisimple linear part. By the argument in the proof of \cite{LWZ}, Theorem 5, it can be written as $R(x)=\Theta(S(x))$, with each entry of $S$ a Laurent monomial. Noting that the rank of $D\Theta$ is necessarily maximal, one sees that the entries of $S$ are first integrals of the normal form. Considering the nontrivial direction of part (c), let $Q_1,\ldots,Q_d$ define independent Laurent monomial first integrals of \eqref{sys}. They span (over $\mathbb Z$) a rank $d$ submodule of $\bar R_\lambda$. Thus for any $Q\in \bar R_\lambda$ there exist integers $\ell\not=0,\ell_1,\ldots,\ell_d$ such that $\ell Q=\sum \ell_iQ_i$. Therefore
    \[
    \left(x^Q\right)^\ell=x^{\ell Q}=\prod_i { x^{\ell_i Q_i}}
    \]
    is a first integral of \eqref{sys}. But then, by general properties of first integrals, so is $x^Q$.

    The proof of part (d) is the same as for Proposition \ref{pro:69}.
\end{proof}

{
\begin{remark}
   Proposition \ref{pro:LWZ} also holds when $d={\rm rank}\,\bar R_\lambda=0$: In this case the eigenvalues $\lambda_1,\ldots,\lambda_n$ are linearly independent over $\mathbb Q$; in particular there are no resonances, and the normal form is linear with diagonal matrix.
\end{remark}
}

For a final observation regarding first integrals, we recall the setting and notation of Bruno's Case 1 (see section \ref{sec:Bruno1}, in particular \eqref{eq:6}). Thus we have
$$
 \lambda=\mu+i \nu,
$$
where $\mu$ and $\nu$ are real vectors, and 
\[
0=\mu_1=\dots =\mu_l< \mu_{l+1}\le \dots \le \mu_n.
\]

Recall from section \ref{sec:Bruno1} that for any $n$-tuple $v=(v_1, \dots, v_n ) $ we denote by 
$v'$ the $l$-tuple $(v_1, \dots, v_l ) $ and by 
$v''$ the $n-l$-tuple $(v_{l+1}, \dots, v_n )$. 
This decomposition has relevance for first integrals:
\begin{lem}\label{FirstintBruno1}
    Given the setting of section \ref{sec:Bruno1} with \eqref{eq:6}, every monomial first integral $x^P$ of $\dot x=A_sx$ satisfies $P''=0$. In other words, every monomial first integral, as well as every formal power series first integral, is constant on the subspace 
    \[
   W:= \left\{\begin{pmatrix}
        0\\ x''
    \end{pmatrix},\quad x''\in \mathbb C^{n-l}\right\}.
    \]
\end{lem}
\begin{proof}
    The condition $\la P,\lambda\ra=0$ implies $\la P,\mu\ra=0$ and thus $\sum_{i>l}p_i\mu_i=0$, whence $p_{l+1}=\cdots=p_n=0$.
\end{proof}
\begin{remark}
    {However, in general $W$ is not the largest subspace on which every formal power series integral is constant; compare Walcher \cite{WalcherNF}, Theorem 1.9 ff.}
\end{remark}

\subsection{Convergence results}
 We now discuss settings in which the existence of sufficiently many (formal or formally meromorphic) first integrals of an analytic system implies the existence of a convergent normalization, by ensuring that some version of Bruno's Condition A holds.
The first result is known 
{(for the case of analytic first integrals see Zhang \cite{XZ}, Theorem 1.1;} for the case of formal first integrals see \cite{LPW}, Theorem 9; for the formally meromorphic case see \cite{LWZ}, Theorem 3), but here we provide a different and simpler proof.
\begin{pro}\label{pro:33}
 Assume that \eqref{sys} is analytic (not necessarily in normal form) and admits $n-1$ independent formally meromorphic first integrals. 
 Then some normal form satisfies Bruno's Simplified Condition A, and there exists a convergent normalizing transformation.
\end{pro}
\begin{proof}
By Proposition \ref{pro:LWZ}, some formal normal form \eqref{nfB1} admits $n-1$ independent Laurent monomial  first integrals $x^{Q_1},\dots, x^{Q_{n-1}}$. 
 Let $I$ be the subspace of $\C^n$ generated by $Q_1,\dots,Q_{n-1}$. 
 Since $Q_1,\dots,Q_{n-1}$ are independent, one has $\dim I=n-1$, and with linear algebra we can write $\C^n=I\oplus \mathbb Cv$, with $v\not=0$, $\la I, v\ra=0$. Since all the $Q_i$ have integer entries, we may take $v\in\mathbb Z^n$. With $\bar 0\in M_\lambda$ and $G_{\bar 0}=\lambda$,
 from \eqref{prPG} we have $\la Q_i, \lambda \ra =0$ for $i=1,\dots, n-1,$ which implies
 \[
  \lambda= cv
 \]
 for some $c\in \mathbb C^*$.
Now take any $P\in M_\lambda$.  By Lemma \ref{lem:it} we have $\la G_P,\, I\ra =0 $, thus $G_P\in\mathbb C v=\mathbb C\lambda$, 
which shows that 
\be \label{Gnm1}
G=(\sum_{P\in M_\lambda} a_P x^P)\lambda\odot x .
\ee

{
To ensure that Bruno's Simplified Condition $A$ is satisfied, it is sufficient to show that if $P\not \in M_\lambda^+$,  then $a_P=0$ for the corresponding $P$ in \eqref{Gnm1}.

Assume that  $P\in M_\lambda$ has a negative entry, say in the $k$-th position.  Then,  there exists an index $i\neq k$ such that
\(
\lambda_i\neq0.
\)
Hence, $(G_P)_i=
a_P\lambda_i$. On the other hand,  since $G_P$ is a coefficient of the normal form,   $ (G_P)_i=0,  $ 
yielding $a_P=0.$ That is, \eqref{Gnm1} is the sum over $P\in M_\lambda^+.$
}

For convergence, there remains to verify Condition $\omega$. But with $\lambda=cv$ we see that $\left<Q,\lambda\right>\in c\mathbb Z$ for any $Q\in {\bf N}$, and therefore $|\left<Q,\lambda\right>|\geq |c|$ whenever $\left<Q,\lambda\right>\not=0$.
\end{proof}
\begin{remark}{\em 
    As known from Francoise \cite{Franc}, maximal formal integrability of an analytic system (with $d<n-1$) does not generally imply the existence of a convergent normalizing transformation. }
\end{remark}
The following extension of Proposition \ref{pro:33} seems to be unavailable in the literature.
{
\begin{pro}\label{pro:33plus}
Let the setting of Bruno's Case 1 (see section \ref{sec:Bruno1}) be given.
Assume that \eqref{sys} is analytic (not necessarily in normal form), and that the eigenvalues $\lambda_1,\ldots,\lambda_l$, $\l\geq 2$ satisfy condition \eqref{eq:6}, with at least {  one } of them nonzero. If the system admits $l-1$ independent formal first integrals, then there exists a convergent normalizing transformation.
\end{pro}
\begin{proof}
    According to \eqref{nfBap} any formal normal form is triangular, in particular one has
    \[
\frac{d}{dt}x'=G'(x').
    \]
    By  Lemma \ref{FirstintBruno1}, the monomial first integrals of the normal form stand in 1-1 correspondence with the monomial first integrals of $G'$. 
    {Moreover, any resonant $P$ is of the form $P=(P',0)$.}  
    As in the proof of Proposition \ref{pro:33} one now finds that necessarily $G'(x')=a(x')\lambda'\odot x'$ with a formal power series $a$, and furthermore one sees that $\lambda_1,\ldots,\lambda_l$ are pairwise commensurable. The assertion now follows from Theorem \ref{Thm:Bruno1*}.
\end{proof}
}

 { The following result, which seems to be unavailable in the literature, extends the convergence criterion from Proposition \ref{pro:33}, which is based on the existence of 
$n-1$ independent formally meromorphic first integrals. Under an additional assumption on the eigenvalues we show that  it is sufficient to assume the existence of only 
$n-2$ independent formally meromorphic first integrals.}

\begin{pro}\label{compint2}
    Assume that \eqref{sys} is analytic (not necessarily in normal form), 
    that $\lambda=\begin{pmatrix}
        \lambda_1\\ \vdots\\ \lambda_n
    \end{pmatrix}$ and $\bar\lambda$ are linearly independent {
    and  the system admits $n-2$ independent formally meromorphic first integrals.  Then any normal form satisfies  Conditions  AS and $\omega$,  and, therefore,}   there exists a convergent normalizing transformation. 
\end{pro}
\begin{proof}
By Proposition \ref{pro:LWZ}, some formal normal form \eqref{nfB1} admits $n-2$ independent Laurent monomial first integrals $x^{Q_1},\dots, x^{Q_{n-2}}$. 
    In this case the solution space of the linear system $\la Q_j, x\ra =0$, $1\leq j\leq n-2$, is spanned by $\lambda$ and $\bar\lambda$. 
 { Now Proposition \ref{nofomerofull} below shows that Condition AS is satisfied.  }
    
    There remains to show Condition $\omega$. Consider the vector space over $\mathbb Q$ that is spanned by the entries of $\lambda$. From $\left<Q_j,\lambda\right>=0$, $1,\leq j\leq n-2$ one sees that this vector space has dimension two; hence there exist complex numbers $\rho_1,\,\rho_2$ and rational vectors $\mu^{(1)},\,\mu^{(2)}$ such that $\lambda=\rho_1\mu^{(1)}+\rho_2\mu^{(2)}$. The normal form conditions and the linear independence of $\lambda,\,\bar\lambda$ remain valid when replacing $\lambda$ by $\rho_1^{-1}\lambda$, hence we may assume that
    \[
    \lambda=\mu^{(1)}+\rho\mu^{(2)},\quad \rho\in \mathbb C\setminus\mathbb R.
    \]
    Since the entries of the $\mu^{(j)}$ are rational, there exists $\delta>0$ (for instance, take the inverse of a common denominator) so that for any $P\in\bf N$ we have $|\left<P,\,\mu^{(j)}\right>|\geq \delta$ whenever $\left<P,\,\mu^{(j)}\right>\not=0$. So we get
    \[
    |\left<P,\,\lambda\right>|\geq  |{\rm Im}(\left<P,\,\lambda\right>)|\geq |{\rm Im}(\rho)|\,\delta\quad\text{whenever  }\left<P,\,\mu^{(2)}\right>\not=0,
    \]
    and 
     \[
    |\left<P,\,\lambda\right>|=  |\left<P,\,\mu^{(1)}\right>)|\geq \delta\quad\text{whenever  }\left<P,\,\mu^{(1)}\right>\not=0\text{  but  }\left<P,\,\mu^{(2)}\right>=0.
    \]
    In particular Bruno's Condition $\omega$ is satisfied. 
    {Thus, by Theorem \ref{BSH} there exists a convergent transformation to normal form. } 
\end{proof}
\begin{exstar}
{\em 
     With $r\geq 2$, $s\geq 2$, let  $p_1,\ldots,p_r,\,q_1,\ldots,q_s$ be nonzero integers, and let $\beta\in \mathbb C\setminus\mathbb R$. Then Proposition \ref{compint2} is applicable to the case when $\lambda=(p_1,\ldots,p_r,\beta q_1,\ldots,\beta q_s)^\top$, since the linear part admits $r+s-2$ independent Laurent monomial first integrals. Moreover, if an analytic system $\dot x=\lambda\odot x+\cdots$ (not necessarily in normal form) admits $r-s-2$ independent formally meromorphic first integrals, then a convergent transformation to normal form exists.}
 
\end{exstar}
\subsection{Case study: Roots of unity}
{In this subsection we consider  system   \eqref{sys} with the matrix 
\be\label{Acycln}
A=\diag(1,\zeta, \zeta^2,\dots, \zeta^{n-1} ),
\ee
where $n>2$ and $\zeta = e^{2\pi i/n}$ is a  primitive $n$-th root of unity. Here Condition $\omega$ holds due to Theorem \ref{prop:omega}.
We consider the formally meromorphic first integrals of such systems.

{By Propositions \ref{formfi}, \ref{pro:69} and \ref{pro:LWZ} we are interested in Laurent monomial first integrals of $\dot x=Ax$, equivalently in integer solutions of the system
\begin{equation}
    \sum_{i=0}^{n-1} m_i \zeta ^i  =0.
\end{equation}
We recall a few facts about roots of unity and cyclotomic field extensions; for details and more see Lang \cite{Lang}, Ch.\ VI, \S 3. 
\begin{itemize}
    \item The minimal polynomial of a primitive $n$th root of unity, $\zeta$, over the rationals $\mathbb Q$ is the $n$th cyclotomic polynomial, $\Phi_n$.
    \item The degree of $\Phi_n$, which is also the degree of the field extension $\mathbb Q(\zeta):\mathbb Q$, is equal to the Euler totient
    \[
    \varphi(n)= \#\left\{ k;\, 1\leq k\leq n,\, {\rm gcd}\,(k,n)=1\right\}.
    \]
    \item When $n$ is a prime number then $\Phi_n(x)=x^{n-1}+\cdots+x+1$.
\end{itemize}
\begin{pro}
    For any  $n>2$, the number of independent formally meromorphic first integrals of $\dot x=Ax$ is equal to $n-\varphi(n)$.
\end{pro}
\begin{proof}
    We abbreviate $d=\varphi(n)$. Let $\zeta$ be a primitive $n$-th root of unity. From $\Phi_n$ one has an integer relation
    \begin{equation}\label{zetamin}
    \zeta^d-\sum_{i=0}^{d-1} m_{i,d}\zeta^i=0
    \end{equation}
    Multiplying repeatedly by powers of $\zeta$ and invoking \eqref{zetamin}, one obtains
    \begin{equation}
        \begin{array}{rcl}
           \zeta^{d+1}& - & \sum_{i=0}^{d-1}m_{i,d+1}\zeta ^i=0 \\
             & \vdots&   \\
             \zeta^{n-1}&-&\sum_{i=0}^{d-1} m_{i,n-1}\zeta ^i=0\\
        \end{array}
    \end{equation}
    with integers $m_{ij}$.
    This system is in echelon form, hence it contains $n-d$ independent solutions. But there cannot exist more than $n-d$ independent solutions; otherwise $1,\zeta,\ldots,\zeta^{d-1}$ would be linearly dependent, whence the minimal polynomial of $\zeta$ would have degree $<d$.
\end{proof}
\begin{cor}
    Let the system \eqref{Asn}, with matrix \eqref{Acycln}, be maximal formally meromorphically integrable. Then Bruno's Condition $A_2$ is satisfied if and only if $n\in \{3,\,4,\,6\}$.
\end{cor}
\begin{proof} 
{The convex hull condition is satisfied for all $n^{\rm th}$ roots of unity, $n\geq 3$.
    Thus the necessary and sufficient condition for the existence of $n-2$ independent integrals is $\varphi(n)=2$, which holds only for these $n$.}
\end{proof}
}
{
\begin{exstar} 
{\em 
    \begin{itemize}
        \item In dimension three, with $\zeta$ a primitive third root of unity, we have $\lambda=(1,\,\zeta,\,\zeta^2)^\top$. In case of maximal integrability, Proposition \ref{pro:33} applies.
        \item In dimension four, let $\lambda=(1,\,i,\,-1,\,-i)^\top$, noting that $i$ is a primitive fourth root of unity. In case of maximal integrability, Proposition \ref{compint2} applies with first integrals $x_1x_3$ and $x_2x_4$.
        \item 
{In dimension six, we have $\zeta^2-\zeta+1=0$ from the sixth cyclotomic polynomial, and furthermore $\zeta^{2+j}-\zeta^{1+j}+\zeta^j=0$ for all $j$ modulo 6. This is a case for Proposition \ref{compint2}: the linear part admits Laurent monomial first integrals
\[
\dfrac{x_1x_3}{x_2},\,\dfrac{x_2x_4}{x_3},\,\dfrac{x_3x_5}{x_4},\,\dfrac{x_4x_6}{x_5},\,\dfrac{x_5x_1}{x_6},\,\dfrac{x_6x_2}{x_1},
\]
with any four of these independent. So here maximal formally meromorphic integrability (as well as Condition $A_2$) holds when four independent formally meromorphic first integrals exist.}
    \end{itemize}
    }
\end{exstar}
}
In prime dimension, maximal formal integrability is equivalent to $\psi=x_1\cdots x_n$ being a first integral, since the $n$-th  cyclotomic polynomial is irreducible over the rationals.
{We give an explicit representation of maximal formally integrable systems in the latter case.
\begin{pro}\label{prop44}
{Let $n$ be a prime number. Then the  system \eqref{nfB1} in normal form
admits the first integral $\psi(x)=x_1\cdots x_n$ if and only if it has a representation \eqref{Gdec}
with
\[
\lambda^{(r)}=(1,\zeta^r,\zeta^{2r},\dots,\zeta^{(n-1)r})\in\mathbb{C}^n,\quad 1\leq r\leq n-1.
\]
}
\end{pro}
\begin{proof}
The system admits the first integral $\psi$ if and only if it has a representation
\[
G(x)=(\lambda^{(1)}\odot x)+(g(x)\odot x),
\]
where $g(x)=(g_1(x),\dots,g_n(x))^\top$ satisfies 
\[
\sum_{j=1}^n g_j(x)\equiv 0,
\qquad g(0)=0.
\]
Letting
\[
V:=\Big\{u\in\mathbb{C}^n:\ \sum_{j=1}^n u_j=0\Big\},
\]
we equivalently have
\[
g(x)=\sum_{Q> 0} U_Q\,x^Q;\quad \text{all  } U_Q\in V.
\]

To see that the $\lambda^{(r)}$, $1\leq r\leq n-1$ form a basis of $V$, consider the cyclic shift operator
\[
S:\,V\to V,\quad \begin{pmatrix} v_0\\ v_1\\ \vdots\\ v_{n-1}\end{pmatrix} \mapsto \begin{pmatrix} v_1\\ \vdots\\ v_{n-1}\\v_0\end{pmatrix}.
\]
A direct computation shows that each $\lambda^{(r)}$ is an eigenvector of $S$ with eigenvalue $\zeta^r$,
and the eigenvalues $\{\zeta^r:\ r=0,1,\dots,n-1\}$ are pairwise distinct.

The assertion follows.
\end{proof}
\subsection{{Further aspects of formally meromorphic integrability}}\label{subsmoremero}
{
In this subsection we first consider formal power series vector fields admitting a representation as linear combination of linear vector fields with formally meromorphic coefficients.
In detail, for given $r<n$, we consider 
\be\label{lamlinco}
G=\begin{pmatrix}
    g_1\\ \vdots \\g_n
\end{pmatrix}=\begin{pmatrix}
    \lambda_{11} x_1\\ \vdots\\ \lambda_{1n}x_n
\end{pmatrix}\, v_1+\cdots+\begin{pmatrix}
    \lambda_{r1} x_1\\ \vdots\\ \lambda_{rn}x_n
\end{pmatrix}\,v_r
\ee
with complex $\lambda_{ij}$ and (a priori) formally meromorphic $v_j$. However, note that we have formal power series $g_i$ on the left hand side, with $g_i(0)=0$. 
{Such representations appear naturally in Proposition \ref{pro:LWZ}. On the other hand,  Conditions A, AS and AL in Propositions \ref{BSH} and \ref{stolopart} about convergence involve such linear combinations but {\em with formal power series coefficients}. 
{Here we intend to show (among other results) that maximal formally meromorphic integrability already implies this additional condition on the coefficients.}}

We set
\[
\Lambda:=\begin{pmatrix}
    \lambda_{11}x_1 &\cdots&\lambda_{r1}x_1\\
    \vdots & & \vdots\\
    \lambda_{1n}x_n&\cdots&\lambda_{rn}x_n
\end{pmatrix}, \quad \Lambda^0:=\begin{pmatrix}
    \lambda_{11} &\cdots&\lambda_{r1}\\
    \vdots & & \vdots\\
    \lambda_{1n}&\cdots&\lambda_{rn}
\end{pmatrix}
\]
and require that $\Lambda^0$ has rank $r$ (thus $\Lambda$ has rank $r$ whenever all $x_i\not=0$).
\begin{remark}\label{rem:req1}
When $r>1$ such representations are not unique; in fact one may replace $\Lambda^0$ by any complex matrix that is equivalent under Gauss column operations over $\mathbb C$.
\end{remark}
{
\begin{pro}\label{cramerprop}
    Given a representation \eqref{lamlinco}, let $1\leq i_1<\cdots<i_r\leq n$ such that the rows $i_1,\ldots,i_r$ of $\Lambda^0$ are linearly independent. 
    \begin{enumerate}[(a)]
        \item Then $v_i\cdot x_i$ is a formal power series for $1\leq i\leq n$.
        \item Moreover $x_j$ divides $g_j$ for every index $j\not\in\{i_1,\ldots,i_r\}$.
    \end{enumerate}
\end{pro}
\begin{proof}
    We may assume without loss of generality that $i_1=1,\ldots,i_r=r$. Let 
    \[
    \bar\Lambda:=\begin{pmatrix}
    \lambda_{11}x_1 &\cdots&\lambda_{r1}x_1\\
    \vdots & & \vdots\\
    \lambda_{1r}x_r&\cdots&\lambda_{rr}x_r
\end{pmatrix}
    \]
    and $\bar\Lambda_k$ the matrix that is obtained from $\bar\Lambda$ by replacing column $\# k$ by $\begin{pmatrix}
        g_1\\ \vdots\\g_r
    \end{pmatrix}$.
    By Cramer's rule one has $v_k=\det\bar\Lambda_k/\det\bar\Lambda$, with 
     \[
    \det\bar\Lambda = x_1\cdots x_r \cdot\det\Lambda^0 \not=0,
    \]
    and 
    \[
    \det\bar\Lambda_k= x_1\cdots x_{k-1}x_{k+1}\cdots x_r\,\Theta_k,
    \]
     with some formal power series $\Theta_k$ (to  verify the latter, expand the determinant along column $\#k$). This shows the first assertion.
    From \eqref{lamlinco} the $j$-th row may be rewritten as
    \[
\det \Lambda^0 x_1\cdots x_r\,g_j=x_j\left( \sum_{k=1}^r \lambda_{kj}x_1\cdots x_{k-1}x_{k+1}\cdots x_r\,\Delta_k\right)
\]
with a formal power series inside the bracket, hence 
   $x_j$ divides $\det\Lambda^0\cdot g_j$, for every $j>r$. 
\end{proof}
}
We call an index $j$ {\em exceptional} (with respect to the representation \eqref{lamlinco}) if the matrix obtained from $\Lambda$ by omitting row $\# j$ has rank $<r$. According to the Proposition, for any non-exceptional index $k$ one has that $g_k=x_k\cdot s_k$ with some formal power series $s_k$. 

{
{
We obtain a clear understanding of exceptional indices by linear algebra. With Gauss column operations one sees the first statement of the following:
{
\begin{lem}\label{LamMulem}
\begin{enumerate}[(a)]
    \item 
There is an invertible matrix $P$ such that (up to renumbering)
\[
\Lambda^0\,P=\begin{pmatrix}
    * & & \cdots& *\\
    \vdots&&&\vdots\\
    * & & \cdots& *\\
     1 &0 & \cdots& 0\\
    0&1&&\vdots\\
    \vdots& &\ddots&0\\
    0& \cdots & 0&1
\end{pmatrix} \in \mathbb R^{n\times r}.
\]
If the corresponding columns of $\Lambda$ commute with some vector field $\mu\odot x$, then the same holds for $\Lambda\,P$. In particular, the normal form property is preserved.
\item An index $j$ is exceptional if and only if $j>n-r$ and column $\# (j-n+r)$ contains only one nonzero entry.
    \end{enumerate}
\end{lem}
\begin{proof}
To prove part (b), first note that omitting a row with number $\leq n-r$ will not decrease the rank, so $j>n-r$. Furthermore, omitting a row with number $j>n-r$ does not decrease the rank if and only if column $\#(j-n+r)$ of the remaining matrix has a nonzero entry.
\end{proof}
}
}}

{
\begin{pro}\label{Wprop}
    Given a vector field $G$ as in \eqref{lamlinco}, let $j_1,\ldots,j_k$ be the non-exceptional indices. Then every subspace defined by $x_{j_\ell}=0$, and in particular the subspace $W$ defined by $x_{j_1}=\cdots=x_{j_k}=0$ is invariant for $\dot x=G(x)$. 
\end{pro}
\begin{proof}
    We may assume that the exceptional indices are $k+1,\ldots,n$ and that $\Lambda$ has the form in Lemma \ref{LamMulem}. Thus we see from the differential equation that 
    \[
    \dot x_i=x_i(\sum_{j=1}^k\lambda_{ji}v_j),\quad 1\leq i\leq k,
    \]
    and the $v_k$ are formal power series by Proposition \ref{cramerprop}. The claim follows with a standard invariance criterion.
\end{proof}
}
{
{
Now we apply these results to normal forms in the setting of Proposition \ref{pro:LWZ}.
{
\begin{pro}\label{nofomerofull}
 Let the formal vector field  $G$ in \eqref{nfB1} be in Poincar\'e-Dulac
  normal form, and satisfy the maximal formally meromorphic integrability condition. Moreover let
 $\bar R_\lambda $ have basis $Q_1,\ldots,Q_{n-r}$ (thus rank $n-r$), and let \eqref{ladec} be a corresponding additive decomposition with $\lambda^{(1)},\ldots,\lambda^{(r)}$ rational isoresonant.
 \begin{enumerate}[(a)]
\item If 
this representation 
admits the exceptional index $j$, then entry $\# j$ of every $Q_i$ is zero. In other words, every Laurent monomial first integral, and every formally meromorphic first integral, of \eqref{nfB1} is independent of $x_j$.
\item $G$ admits a representation 
\[
    G(x)=\sum_{i=1}^{r} (\gamma_i+\check s_i(x))\cdot \lambda^{(i)}\odot x
    \]
  with {\em formal power series} coefficients
  \[
    \check s_i(x)=\sum_{P\in M_\lambda^+ }  \sigma_{i,P} x^P, \quad   s_i(0)=0, \quad i=1,\dots, r.  
\]
 \end{enumerate}
\end{pro}
\begin{proof}
    We may assume that $j\in\{k+1,\ldots,n\}$ and $\Lambda$ as in Lemma \ref{LamMulem}. The $Q_i$ (seen as row vectors) lie in the left kernel of $M$, hence must have entry $\# j$ equal to zero. This shows part (a). For part (b), one has a representation as stated, but a priori the representation of $\check s_i$ involves $P \in M_\lambda$. One needs to show that it suffices to restrict to summation over nonnegative $P$. For entry $\#j$ of $G$ the assertion follows when $j$ is non-exceptional, by Proposition \ref{cramerprop}. If $j$ is exceptional, then Lemma \ref{LamMulem} shows that $\lambda^{(j-n+r)}=e_j$; the canonical basis vector, and entry $j$ of the differential equation reads as 
    \[
    \dot x_j=\sum_R \alpha_R x^R=\sum_R\alpha_R\dfrac{x^R}{x_j}\cdot x_j; \quad R\in \mathbb Z_{\geq 0}^n,\,\alpha_R\in \mathbb C.
    \]
    Whenever $\alpha_R\not=0$, then $x^R/x_j$ is a Laurent monomial first integral of the linear part $\lambda\odot x$, and thus of the normal form $G$ by maximal formally meromorphic integrability. The assumption that $x_j$ does not divide $x^R$ implies that there is a Laurent monomial first integral that depends on $x_j$; a contradiction to (a). We see that $x_j$ divides $x^R$ whenever $\alpha_R\not=0$.
    
\end{proof}
}}
\begin{remark}
    In the setting of Proposition \ref{nofomerofull}, the subspace $W$ from Proposition \ref{Wprop} is the largest subspace on which every formally meromorphic first integral of the normal form is constant. In some cases this subspace corresponds to the decomposition from Remark \ref{rem:echelon}, but this is not true in general.
\end{remark}

A variant of the proof of Proposition \ref{nofomerofull} yields:
\begin{pro}\label{nofodiohull}
 Let the formal vector field  $G$ in \eqref{nfB1} be in Poincar\'e-Dulac
  normal form, and let $\lambda$ admit an additive decomposition with $\lambda^{(1)},\ldots,\lambda^{(r)}$ satisfy the diophantine hull condition.
Then $G$ has a representation 
\[
    G(x)=\sum_{i=1}^{r} (\gamma_i+\check s_i(x))\cdot \lambda^{(i)}\odot x
    \]
  with {formal power series} coefficients
  \[
    \check s_i(x)=\sum_{P\in M_\lambda^+ }  \sigma_{i,P} x^P, \quad   s_i(0)=0, \quad i=1,\dots, r.  
\]
\end{pro}
\begin{proof} We follow the pattern in the proof of Proposition \ref{nofomerofull} (b). Thus let $j$ be an exceptional index, then Lemma \ref{LamMulem} shows that $\lambda^{(j-n+r)}=e_j$; the canonical basis vector, and entry $j$ of the differential equation reads as 
    \[
    \dot x_j=\sum_R \alpha_R x^R=\sum_R\alpha_R\dfrac{x^R}{x_j}\cdot x_j; \quad R\in \mathbb Z_{\geq 0}^n,\,\alpha_R\in \mathbb C.
    \]
    Now, if $\alpha_R\not=0$, then $\la R-e_j,\lambda\ra=0$ since $G$ is in normal form, and the diophantine hull condition implies $\la R-e_j,\lambda^{(j-n+r)}\ra=0$. With $R=\sum m_ke_k$, $m_k\in\mathbb Z_{\geq 0}$, this implies $\la (m_j-1)e_j,e_j\ra=0$, hence $m_j=1$, and $x_j$ divides $x^R$.
\end{proof}
\begin{remark}\label{rem:A2revisit}
This result is applicable to Bruno's Condition $A_2$. Bruno implies (with no further comment) in \cite{BruB}, p.194 that a normal form satisfying Condition $A_2$ has a representation
\[G(x)=a(x)\lambda\odot x+b(x)\bar\lambda\odot x\]
with formal power series $a$, $b$. The Proposition provides an explanation, and shows a bit more:
Condition $A_2$, as given by Bruno, consists of two requirements. The first is that $\lambda$ and $\bar\lambda$ are linearly independent (note that the diophantine hull condition for $\lambda$ and $\bar\lambda$ is satisfied). The second requirement is that $0$ lies in the interior of the convex hull of $\lambda_1,\ldots,\lambda_n$; we nowsee that it is not necessary.
  \end{remark}}
  }

We now turn to another aspect of maximal formally meromorphic integrability, with regard to convergence criteria for normalization of an analytic vector field \eqref{Asn}. The relevance of nontrivial commuting analytic vector fields (``infinitesimal symmetries'') for convergence was noted in Bruno and Walcher \cite{BrW} (proving convergence in dimension two), and in the expository paper \cite{CiW} by Cicogna and Walcher (exhibiting some classes of higher dimensional systems). N.T.\ Zung  combined commuting vector fields and first integrals to obtain sufficient conditions for convergence, as follows.
 \begin{teo} {\em (N.T.\ Zung \cite{Zung,Zungrat})}\\
     Given an $n$-dimensional analytic system \eqref{Asn}, assume there exist $1\leq p\leq n$ and 
     \begin{enumerate}[(i)]
         \item independent local meromorphic vector fields $F_1=F,\,F_2,\ldots,F_p$ that commute pairwise; and
         \item independent local meromorphic common first integrals $h_1,\ldots,h_{n-p}$ of the $F_i$.
     \end{enumerate}
     Then there exists a convergent transformation of $\eqref{Asn}$ to normal form.
 \end{teo}
 \begin{remark}\label{rem:ranknul}
 {
 \begin{enumerate}[(a)]
 \item We briefly call (i) and (ii) the Zung conditions, and we speak of the formal Zung conditions when the $F_i$ and $h_j$ are assumed to be only formally meromorphic.
     \item The Zung conditions and the formal Zung conditions describe a  maximality property: Generally, given $n-p$ independent functions, there are at most $p$ independent vector fields which have these as common first integrals; conversely, $p$ independent vector fields have at most $n-p$ common first integrals. Both statements follow from writing out the conditions and recalling the rank-nullity theorem from linear algebra. 
     \item In the case $p=0$, with $n$ pairwise commuting analytic vector fields and (in particular) no resonances among the eigenvalues, a relatively simple proof of Zung's theorem may be obtained, for instance, from Theorem 2.4 and Corollary 2.5 in Cicogna and Walcher \cite{CiW}.
 \end{enumerate}
     }
 \end{remark}
 {
 {We now look at Zung's conditions in the formal framework, to understand where possible formal obstructions to convergence lie.} Here we have:
 \begin{pro}\label{formzung}
     Let an $n$-dimensional formal system \eqref{Asn} be given. Then the formal Zung conditions will hold if and only if the system is maximal formally meromorphic integrable.
 \end{pro}
 \begin{proof}
 {We may assume that the system is in normal form. The ``if'' part is a direct consequence of Proposition \ref{pro:LWZ} (set $p=d$), with the commuting vector fields $F_j(x)=\lambda^{(j)}\odot x$. For the ``only if'' part, given (i) and (ii) in the formal setting, part (a) of Remark \ref{rem:ranknul} shows that there cannot be more than $n-p$ independent formally meromorphic first integrals.}
 \end{proof}
 Thus the only formal obstruction lies in the existence of sufficiently many formally meromorphic first integrals.\footnote{The formal Zung conditions generally do not suffice for convergence. Even in dimension two there always exist nontrivial commuting formal vector fields for $F$ by the normal form property, but there are cases when no convergent normalizing transformation exists.}}
{
\begin{remark}
Proposition \ref{formzung} provides a ``recipe'' to construct systems in normal form that satisfy the formal Zung conditions: Start with an independent set $Q_1,\ldots,Q_{n-p}$ in $\mathbb Z^n$, let $\lambda^{(1)},\ldots, \lambda^{(p)} \in \mathbb Q^n$ be linearly independent such that all $\la Q_i,\lambda^{(j)}\ra=0$, and finally let $\gamma_1,\ldots,\gamma_p\in \mathbb C$  be linearly independent over the rationals.  Now use the procedure in the proof of Proposition \ref{pro:LWZ} (setting $\sum \gamma_k \lambda^{(k)}=\lambda$). Thus one prescribes the first integrals and $\bar R_\lambda$. The choice of the $\gamma_j$ ensures that every first integral of $\lambda\odot x$ must be a common first integral of the $\lambda^{(i)}\odot x$. One arrives at
\[
G(x)=\sum_{j=1}^p(\gamma_j+\check s_j(x))\lambda^{(j)}\odot x;\quad \check s_j(x)=\sum_{Q\in\sum\mathbb Z Q_i}\sigma_{jQ}x^Q.
\]
 This observation also provides a new perspective on the constructions by Yamanaka \cite{Yama}. 
 \end{remark}
 As a byproduct of the above considerations we see that the structure of normal form vector fields is rather transparent in the formally meromorphic setting. Matters are more complicated when one considers formal power series vector fields and first integrals. For some results see Kruff et al.\ \cite{KWZ}; in particular section 4.
 }
\\

Finally,
we return to the scenario of section \ref{sec:stolo} and consider again the simultaneous normalization of pairwise commuting analytic or formal vector fields $F^{(1)},\ldots,F^{(s)}$ on $\mathbb C^n$, given as in \eqref{eq:abellie}, thus 
\begin{equation*}
F^{(k)}(x)=A^{(k)}x+\text{t.h.o} =\lambda^{(k)}\odot x +\text{t.h.o},
\end{equation*}
and linearly independent $\lambda^{(1)},\ldots,\lambda^{(s)}$.\footnote{ Note: This implies that the $F^{(k)}$ are independent.}

From Corollary \ref{stoloform} it is known that there exists a formal transformation which sends the $F^{(k)}$ to commuting vector fields $G^{(k)}$ in simultaneous normal form.

For analytic $F^{(i)}$, Theorem \ref{stolopart} shows convergence when conditions AL and $\omega^{\#}$ are satisfied. We will now show that the existence of independent formally meromorphic first integrals $h_1,\ldots, h_{n-s}$ will force condition AL.

{
{
To proceed, we need the following definition and an auxiliary statement. 

\begin{defin}\label{def:max-family}
We say that the family  \eqref{eq:abellie} is maximal formally meromorphic
integrable, if it admits \(n-s\) independent common formally meromorphic
first integrals.
\end{defin}
By Remark~\ref{rem:ranknul}, a family of \(s\) independent vector
fields on \(\mathbb C^n\) admits at most \(n-s\) independent common
first integrals. Thus Definition~\ref{def:max-family} describes precisely
the case in which this upper bound is attained.

Since the $F^{(k)}$ admit $n-s$ common formally meromorphic first integrals, the normalized vector fields $G^{(k)}$ admit $n-s$ common Laurent monomial first integrals $x^{Q_j}$ with independent $Q_1,\ldots,Q_{n-s}\in \mathbb Z^n$; see Proposition \ref{pro:LWZ}. In particular,
all $\la Q^{(j)},\lambda^{(k)}\ra=0$, and therefore the subspace of all $v\in\mathbb C^n$ with $\la Q_1,v\ra=\cdots=\la Q_{n-s},v\ra=0$ is spanned by $\lambda^{(1)},\ldots,\lambda^{(s)}$.
\begin{lem}
    There exist $\beta_{ij}\in\mathbb C$, $1\leq i,j\leq s$, such that the 
    \[
    \widetilde\lambda^{(i)}:=\sum_{j=1}^s\beta_{ij}\lambda^{(j)},\quad 1\leq i\leq s
    \]
    are linearly independent, and the system $\lambda^{(1)},\ldots,\lambda^{(s)}$ satisfies the rational isoresonance condition with respect to every $\widetilde\lambda^{(k)}$.
    Moreover,
    the system $\widetilde \lambda^{(1)},\ldots,\widetilde\lambda^{(s)}$ satisfies the rational isoresonance condition with respect to every $\widetilde\lambda^{(k)}$.
\end{lem}
\begin{proof}
    Let $\lambda^{(j)}_1,\ldots,\lambda^{(j)}_s\in\mathbb C$ denote the entries of $\lambda^{(j)}$, and set
    \[
    K:=\mathbb Q\left(\left\{ \lambda^{(j)}_\ell;\quad 1\leq j\leq s,\,1\leq \ell\leq n\right\}\right).
    \]
    Since $\mathbb C$ has infinite transcendence degree over $\mathbb Q$, the same is true for $\mathbb C$ over $K$. In particular, the dimension of $\mathbb C$ as a $K$-vector space is infinite. 

    Now let $\gamma_1,\ldots,\gamma_s$ be linearly independent over $K$, and $\mu:=\sum \gamma_j\lambda^{(j)}$. For any $P \in \mathbb Z^n$ one has that all $\la P,\lambda^{(j)}\ra \in K$, and thus
    \[
    0=\la P,\mu\ra=\sum \gamma_j\la P,\lambda^{(j)}\ra\Longleftrightarrow \text{  all  }\la P,\lambda^{(j)}\ra=0.
    \]
    Thus the system of the $\lambda^{(j)}$ satisfies the rational isoresonance condition with respect to $\mu$.

    Choosing the $\beta_{ij}$ linearly independent over $K$, and so that the matrix $(\beta_{ij})$ is invertible, we see the first statement of the lemma. The second statement is then obvious.
\end{proof}
\begin{pro} \label{pro:Zung-form}
    Assume that the system $F^{(1)},\ldots, F^{(s)}$ of commuting formal vector fields in \eqref{eq:abellie}, with simultaneous normal forms $G^{(1)},\ldots, G^{(s)}$, is maximal formally meromorphic integrable. Then condition AL is satisfied.
\end{pro}
\begin{proof}
One may replace the system of the $F^{(k)}$ by
\[
\widetilde F^{(\ell)}=\sum_k \beta_{k\ell}F^{(k)}=(\sum_k\beta_{k\ell}\lambda^{(k)})\odot x+\cdots, \quad 1\leq \ell\leq s
\]
with any invertible matrix $(\beta_{k\ell})$ over $\mathbb C$. By suitable choice of the $\beta_{k\ell}$ we may assume that the $\widetilde\lambda^{(\ell)} =\sum_k\beta_{k\ell}\lambda^{(k)} $ satisfy the conclusion of the lemma above. For the normal forms we have 
\[
\widetilde G^{(\ell)}=\sum\widetilde\lambda^{(\ell)}\odot x+\sum_{j=1}^s a_{\ell j}\widetilde\lambda^{(j)}\odot x
\]
from Proposition \ref{cramerprop}, with formally meromorphic $a_{\ell j}$, and the rational isoresonance condition holds for the $\widetilde\lambda^{(\ell)}$.
Now the assertion follows directly with Proposition \ref{nofomerofull}.
\end{proof}
As a consequence we get the following variant of Zung's convergence theorem:
\begin{cor}
     When Zung's condition (i) is satisfied, but condition (ii) holds only with formally meromorphic first integrals, and Condition $\omega^\#$ holds in addition, then there exists a convergent normalizing transformation.
\end{cor}
 
 In the settings of Propositions \ref{pro:33} and \ref{compint2}, Condition $\omega^\#$ is automatically satisfied in presence of the maximal formal integrability condition. But beyond, there seem to be no more such cases.
}}

\bigskip
\bigskip

{\bf Acknowledgements}. This paper was initiated during a research stay of SW at the Center for Applied Mathematics and Theoretical Physics, University of Maribor. SW thanks the Center for its hospitality and for creating a congenial research atmosphere.
VR acknowledges the support  by the Slovenian Research  and Innovation  Agency (core research program P1-0306)
and 	by the project  101183111-DSYREKI-HORIZON-MSCA-2023-SE-01 “Dynamical Systems and Reaction Kinetics Networks”.

\end{document}